\documentclass[ejsv2,preprint,noshowframe]{imsart}
\RequirePackage[authoryear]{natbib}
\RequirePackage[colorlinks,citecolor=blue,urlcolor=blue,linkcolor=blue]{hyperref}
\RequirePackage{graphicx,booktabs,array,enumitem,mathtools,amssymb,microtype}
\startlocaldefs
\theoremstyle{plain}
\newtheorem{theorem}{Theorem}[section]
\newtheorem{proposition}[theorem]{Proposition}
\newtheorem{lemma}[theorem]{Lemma}
\newtheorem{corollary}[theorem]{Corollary}
\theoremstyle{definition}

\newcommand{\E}{\mathbb E}\newcommand{\Pp}{\mathbb P}
\newcommand{\kl}{\operatorname{kl}}

\endlocaldefs
\hypersetup{pdftitle={Sequential Certification of Partially Identified Product Targets: Information Limits and Adaptive Sampling},pdfauthor={Abdallah Skiredj; Abderrahman Skiredj; Amal El Fallah Seghrouchni}}
\begin{document}
\begin{frontmatter}
\title{Sequential Certification of Partially Identified Product Targets: Information Limits and Adaptive Sampling}
\pdfsubject{Sequential certification and partial identification}
\runtitle{Sequential certification of product targets}
\begin{aug}
\author[A]{\fnms{Abdallah}~\snm{Skiredj}\ead[label=e1]{abdallah.skiredj@um6p.ma}},
\author[B]{\fnms{Abderrahman}~\snm{Skiredj}\ead[label=e2]{abderrahman.skiredj@um6p.ma}}
\and
\author[A]{\fnms{Amal}~\snm{El Fallah Seghrouchni}\ead[label=e3]{amal.elfallah-seghrouchni@um6p.ma}}
\address[A]{AI Movement, Mohammed VI Polytechnic University (UM6P), Rabat, Morocco\printead[presep={,\ }]{e1,e3}}
\address[B]{College of Computing, Mohammed VI Polytechnic University (UM6P), Benguerir, Morocco; Ministry of AI, Digital Transition and Administrative Reform, Rabat, Morocco\printead[presep={,\ }]{e2}}
\runauthor{A. Skiredj et al.}
\end{aug}
\begin{abstract}
We study the cost of certifying a scalar target driven by a product of centered Bernoulli
means when an unobserved terminal gate leaves irreducible uncertainty. Each
factor can be queried independently at a price. We derive the exact identified
interval and determine when a requested accuracy is feasible. At every fixed
strictly feasible interior instance, the optimal high-confidence expected cost
is $(1+o(1))\log(1/\delta)/D^*$, where $D^*$ is an answer-and-allocation
information game; a globally honest procedure attains this benchmark. At
product cancellation, we solve the game through a one-dimensional inverse-KL
equation and explicit cost weights. For two tickets, quartic projections and
rational bounds give checkable global optimality gaps. Joint
confidence--precision limits expose a failure of the analyzed fixed pilot: it
can cost asymptotically at least twice the optimum. A charged, precision-aware
pilot with checkpoint certification restores optimal expected cost in exact
arithmetic when
$\log\{e+\log(1/g)\}=o(\log(1/\delta))$ for effective precision $g$,
including exponential and certain double-exponential precision scales. Pooled rescue
preserves validity and cannot increase observation cost on coupled paths.
Implemented auditors verify returns by rational inequalities; reproducible
experiments show both savings and finite-confidence overhead. These guarantees
concern the specified independent-reference model and its scalar target.
\end{abstract}
\end{frontmatter}

\section{Introduction}
\label{sec:introduction}
Repeatedly running a simulation does not by itself validate the rules that
convert its actions into outcomes. We study the number of costly reference
observations needed to certify one scalar conclusion of such a simulation.
An adjudicator here is an outcome-generating rule; the reference is a
specified statistical standard. Our motivating contrast is
$\Delta=p_{11}-p_{10}-p_{01}+p_{00}$, where $p_{ij}$ is a terminal-event
probability under policy $i$ and adjudicator $j$. Taking the inexpensive
contribution $C_{\rm cheap}=p_{11}-p_{01}$ as known leaves
$F=p_{00}-p_{10}$ to certify.

The mathematically complete model in this paper consists of independent
Bernoulli reference tickets and an unobserved terminal gate. It captures an
interaction between observation, subsequent action, and terminal outcome, while
remaining small enough to obtain explicit information formulas. It does not
resolve arbitrary partially observed games or optimize live human intervention.

Two kinds of uncertainty must be distinguished. Sampling uncertainty can be
reduced by buying more reference observations. The unobserved gate, however,
leaves several target values compatible with exactly the same observations,
even with unlimited data. A certificate must cover this entire identified
interval. The first question is therefore whether the requested accuracy is
possible; only then is it meaningful to optimize its observation cost.
The second difficulty is cancellation: a product can vanish while its
individual factors remain unknown. An allocation learned adequately at a
fixed precision can become inefficient when the requested precision also
shrinks. This is the phenomenon addressed by the joint-limit theorem.

The contribution has three parts. First, the partially identified product
admits an explicit cancellation information game and allocation, within a
characteristic-cost formulation that includes infeasible parameters in its
honesty domain. Second, uniform product-root stability supports a
precision-aware pilot and a joint expected-cost theorem under an
iterated-logarithm condition, with every pilot observation charged.
Third, two-ticket projections admit computable rational global certificates,
which connect the theory to an executable auditor.
The checkpoint, pilot-localized and pooled certification propositions are
specializations or combinations of established concentration principles;
their role is to support the main theorem and the implementation.
We also separate the price of expert responses from the price of finding
auditable histories under a restricted discovery protocol.

\begin{figure}[t]
\centering
\includegraphics[width=\linewidth]{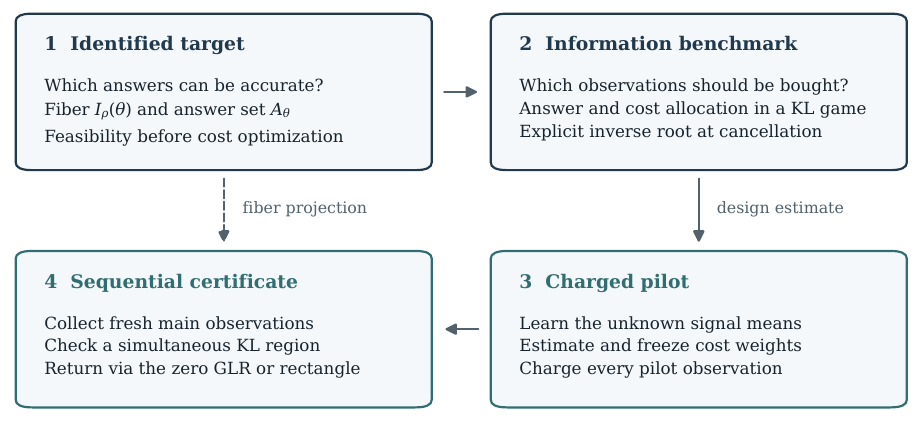}
\caption{Logical structure of the certification problem. The observable model
determines the identified interval and acceptable answers. The information
game supplies a cost benchmark and a target allocation, which the charged
pilot estimates. Fresh main observations then justify stopping through a
simultaneous confidence region. The general game applies at fixed strict
instances; the explicit joint-efficiency theorem concerns cancellation arrays.
The diagram represents mathematical dependencies, not a claim that the
unknown true parameter is available to the algorithm.}
\label{fig:overview}
\end{figure}

\paragraph{How to read the paper.}
Section~\ref{sec:productmodel} introduces the model and a worked example.
Section~\ref{sec:generalnew} defines the information benchmark;
Sections~\ref{sec:explicitnew}--\ref{sec:quarticnew} make it explicit and
computable. Section~\ref{sec:mixturenew} supplies the statistical certificate.
Sections~\ref{sec:jointnew}--\ref{sec:localized} explain why learning the
allocation matters and culminate in Theorem~\ref{thm:loglogjoint}.
Section~\ref{sec:experimentsnew} tests the executable two-ticket procedures.
Figure~\ref{fig:overview} connects these steps. Readers interested primarily
in the main asymptotic result can follow the model, the inverse-KL formula,
and the precision-aware policy before turning to the complete proofs.

\paragraph{Closest prior work.}
Answer-dependent alternatives are central to \citet{degenne2019}.
\citet{poiani2026} already establish lower bounds and optimal algorithms for
regular pure-exploration problems with infinitely many answers. Our general
cost theorem is a specialization and extension to possibly empty acceptable
sets, with honesty required even where termination is impossible; it is not a
claim to invent continuous-answer optimality.
Products of means and cost-aware pilot allocations have substantial precedents
\citep{hardwick1992,zheng1998,benkamra2012,benkamra2015}, including
allocation with a zero factor \citep{benkamra2022}; degeneracy and higher-order inference are also
established \citep{chenfang2019}. Mixture confidence regions are classical
\citep{kaufmannkoolen2021,howard2021}. Recent work by \citet{wang2026logs}
studies optimal certification of policy improvement from logged episodes;
our experiment instead permits individually priced reference tickets and
returns a numerical answer valid for an entire nuisance interval.
Published finite-confidence analysis of Track-and-Stop is also available
\citep{poiani2026finite}; its finite-answer formulation and instance-dependent
remainder do not directly give the moving-product guarantee proved here.
\citet{benkamra2015}, equations (1), (11)--(12) and Theorems 1--3,
study the allocation-variance criterion for a product of positive Bernoulli
proportions at a fixed total sample size. They derive the proportional
allocation, a second-order variance benchmark, and nearly second-order
three-stage rates, including a probability remainder of order
$o_P(\log\log T/T^2)$. Pilot limitations, multistage repair and ordinary
variance-optimal product allocation are therefore precedents, not contributions
of this paper. Our regular-point limit below recovers the corresponding
variance allocation. Their stated results fix interior positive proportions;
they do not supply our answer-dependent KL game, nuisance-fiber honesty,
or uniform expected stopping-cost guarantee on moving cancellation arrays.
Here a centered factor can vanish at an interior Bernoulli parameter.
The distinction is in these losses, domains and quantifiers, not in the
use of products or pilot observations.

Related policy-evaluation and identification formulations
\citep{russo2025,ariu2026,duarte2024,raghavan2026} use different observation
protocols or targets. A comparison of their algorithms under a common
protocol is outside the present study. The contribution claimed here is
the explicit product geometry and its uniform sequential cost analysis.

\section{Experiment, identified interval, and feasibility}
\label{sec:productmodel}
\subsection{Observable model and irreducible uncertainty}
A query to ticket $j\in\{1,\ldots,r\}$ costs $c_j>0$ and returns a fresh
$\operatorname{Bern}(\theta_j)$ observation. Queries may be adaptive;
$\theta\in(0,1)^r$ in efficiency statements. Define
\[
 P(\theta)=\prod_{j=1}^r(2\theta_j-1),\qquad 0\le\rho<\tfrac12.
\]
There are independent signs with these means. One policy chooses their parity,
and a second its negative; an event occurs when the chosen sign is positive
and an unobserved gate succeeds. The gate has a common mean
$\eta\in[\rho,1-\rho]$, with history-dependent deviations in $[-\rho,\rho]$.
All these means remain in $[0,1]$. Averaging deviations over the two parity
classes gives exactly
\[
 F=\eta P+e,\qquad \eta\in[\rho,1-\rho],\quad e\in[-\rho,\rho].
\]
The gate is never queried. Thus all its values induce the same data law.

\begin{proposition}[Exact observable fiber]\label{thm:robustfiber}
The compatible targets at $\theta$ are exactly
\[
 I_\rho(\theta)=[L_\rho(P),U_\rho(P)],\qquad R_\rho(t)=\rho+(1/2-\rho)|t|.
\]
\[
 L_\rho(t)=t/2-R_\rho(t),\qquad U_\rho(t)=t/2+R_\rho(t).
\]
Both endpoints are nondecreasing. The outputs accurate within $\varepsilon$
for every compatible gate form
\[
 A_\theta=[U_\rho(P)-\varepsilon,L_\rho(P)+\varepsilon].
\]
This interval has positive length exactly when $R_\rho(P)<\varepsilon$.
\end{proposition}
\begin{proof}
For fixed $\eta$, the deviation difference fills $[-\rho,\rho]$: each
policy-selected parity class may use a constant deviation of the appropriate
sign. Explicitly, for any $e\in[-\rho,\rho]$, use deviations $e$ and $-e$
on the positive and negative parity classes. Their weighted difference is
$[(1+P)/2]e-[(1-P)/2](-e)=e$; both gate means $\eta\pm e$ lie in $[0,1]$. The union of $[\eta P-\rho,\eta P+\rho]$ over the allowed $\eta$
has the displayed endpoints. Their slopes are $\rho$ and $1-\rho$ on the two
sides of zero. Intersecting the accuracy intervals
$[f-\varepsilon,f+\varepsilon]$ over the fiber gives $A_\theta$.
\end{proof}

An auditor is \emph{globally honest} if, for every signal vector and every
compatible gate,
$\Pp(\tau<\infty,|\widehat F-F|>\varepsilon)\le\delta$.
If $N_j(t)$ counts queries to ticket $j$ by time $t$, total observation
expenditure at return is $C_\tau=\sum_j c_jN_j(\tau)$.
Nontermination has infinite cost. Let $B^*(\theta;\varepsilon,\delta)$ be the
infimum of expected total cost over honest procedures with finite expected
cost at every strictly feasible interior instance. Throughout the paper
$\rho<\varepsilon<1$ and $0<\delta<1/4$, unless stated otherwise.
Write $\ell=\log(1/\delta)$; the joint-limit sections use $L$ for this
log-confidence scale.

\subsection{A worked example and statistical vocabulary}
An \emph{identified interval} contains every target compatible with a fixed
observable law. A \emph{confidence region} instead describes uncertainty
about that observable law from finitely many samples. The auditor projects
the latter through the former. Thus the same return must remain accurate
for both sampling uncertainty and every unobserved gate still compatible
with the data; the gate itself is never estimated from these tickets.

For example, take two tickets, $\rho=.05$ and $\varepsilon=.18$.
At $\theta=(.5,.8)$, the centered product is zero, so
$I_\rho=[-.05,.05]$ and $A_\theta=[-.13,.13]$.
Zero is an acceptable answer, although a finite-data certificate must still
exclude nearby observable alternatives where it would fail.
At $\theta=(.7,.8)$, the product is $.24$, giving
$I_\rho=[-.038,.278]$ and $A_\theta=[.098,.142]$.
Certification remains feasible, but zero is no longer acceptable.
This distinction explains why the general benchmark optimizes the returned
answer as well as the allocation. The zero-answer threshold
$g=(\varepsilon-\rho)/(1-\rho)=13/95$ is stricter than feasibility alone.

We use \emph{cancellation} for $P(\theta)=0$, equivalently at least one
ticket with $\theta_j=1/2$; such a ticket is called \emph{fair}.
The number of fair coordinates determines the local order of the product.
The following boundary result shows why a strict feasibility margin is
needed even though the target is scalar.

\begin{proposition}[The boundary cannot be reliably recognized]\label{prop:boundarynew}
At an interior point with $R_\rho(P)\ge\varepsilon$, every honest procedure
satisfies $\Pp(\tau<\infty)\le2\delta$.
At every strict point $R_\rho(P)<\varepsilon$, honest finite-mean termination
is possible.
\end{proposition}
\begin{proof}
For $R>\varepsilon$, every output fails for at least one of the two
indistinguishable endpoint targets; sum the error probabilities. For equality,
the only common acceptable output is $a=P/2$. Outputs below and above $a$
already fail at opposite endpoints. To include an output exactly at $a$, move
an endpoint outward by an arbitrarily small observable perturbation. If $P>0$,
increase $P$ for the upper endpoint; if $P<0$, decrease it for the lower
endpoint. Such perturbations exist at every interior nonzero product. If
$P=0$, arbitrarily small products of either sign exist, and the extreme targets
are $\pm[\rho+(1-\rho)|P|]$. Split outputs at $a$, assigning equality to the
half excluded by the perturbed endpoint. For example, when $P>0$, the outward upper endpoint excludes
$\{\tau\le m,\widehat F\le a\}$; the base lower endpoint excludes
$\{\tau\le m,\widehat F>a\}$. For $P<0$, reverse these inequalities.
For a bounded number of queries,
transcript probabilities are continuous in the Bernoulli means, including under
adaptive query selection. Let the perturbation vanish, then let the bound on
the query count increase. Each half has probability at most $\delta$.
For strict points, sample every ticket and project a simultaneous shrinking
confidence rectangle through the fiber formula. A positive width margin and
exponential concentration imply a finite expected stopping cost.
\end{proof}

\subsection{Interpretation and scope of the guarantees}
Certifying $F$ is not a certificate against every form of trajectory drift.
With one fair signal and exact sharing, $F=0$ for every gate mean, while
$p_{00}=p_{10}=\eta/2$. Gates $\eta=0$ and $\eta=1$ have identical accessible
tickets and different event probabilities, separated by $1/2$. Their full
trajectory laws have total variation distance at least $1/2$. Therefore no
bound on this contrast alone guarantees that simulated perceptions stay close
to reality or to every reference trajectory. A trajectory-level intervention
question requires an additional metric and a controlled dynamical model.

\paragraph{Assumption and implementation map.}
Theorems~\ref{thm:fullgame} and \ref{prop:plugin} fix an interior strict instance
while confidence increases. Theorems~\ref{thm:jointnew}, \ref{thm:adaptivejoint} and
\ref{thm:loglogjoint} instead require exact cancellation, a fixed compact
interior cube, fixed dimension and costs, and their stated array conditions.
All validity statements assume fresh independent stationary reference responses
and known costs and sharing radius. The general catalog is an existence
construction; the executable learned and adaptive-pilot policies and quartic
validator have two tickets. The reused fixed pilot is supported by
Corollary~\ref{cor:reuse}; the random adaptive pilot uses fresh main data.
For an application, uncertainty in the inexpensive contrast contribution must
also be budgeted: certify it within $\varepsilon_C$ at risk $\delta_C$ and
certify the reference component within $\varepsilon_F$ at risk $\delta_F$.
The triangle inequality and a union bound give contrast error at most
$\varepsilon_C+\varepsilon_F$ with risk at most $\delta_C+\delta_F$,
without requiring independence between the two certificates.

\section{The exact finite-precision characteristic cost}
\label{sec:generalnew}
\subsection{Information per unit of expenditure}
For $p,q\in(0,1)$, Bernoulli relative entropy is
\[
 \kl(p,q)=p\log\frac pq+(1-p)\log\frac{1-p}{1-q}.
\]
Natural logarithms are used throughout. Endpoint values are interpreted by
limits, with $0\log0=0$. Relative entropy measures the mean log-likelihood
evidence supplied by one observation for distinguishing $p$ from $q$.
Because a query to ticket $j$ costs $c_j$, a cost fraction $w_j$ buys
approximately $Cw_j/c_j$ observations out of expenditure $C$.
The weighted KL sum below is therefore information \emph{per unit cost}.

The hardest alternative is one that makes the chosen answer inaccurate
while producing little statistical evidence against itself. A good design
makes even this alternative expensive to confuse with the true law.
Allowing several accurate answers adds an outer choice: the easiest answer
to certify need not be the midpoint of the identified interval. The true
$\theta$ appears in this benchmark for analysis; implementable procedures
must learn their design from observations.

For a possible numerical answer $a$, define its \emph{strict bad set}
\[
 \mathcal B(a)=\{x\in(0,1)^r:L_\rho(P(x))<a-\varepsilon
                 \text{ or }U_\rho(P(x))>a+\varepsilon\}.
\]
For cost fractions $w\in\Delta_r$, put
\begin{equation}\label{eq:generalanswergame}
 d(\theta,w,a)=\inf_{x\in\mathcal B(a)}
        \sum_j\frac{w_j}{c_j}\kl(\theta_j,x_j),\qquad
 D^*(\theta)=\sup_{a\in A_\theta}\sup_{w\in\Delta_r}d(\theta,w,a).
\end{equation}
The infimum uses all interior observable alternatives, including those at
which no common accurate answer exists. Removing those alternatives would
weaken the required global honesty. The value $D^*$ is finite and positive
at every strict point. Indeed a strictly acceptable answer and positive weights
have positive distance from their bad set. Conversely an interior alternative
with product close enough to $+1$ or $-1$ excludes any answer in
$[-\varepsilon,\varepsilon]$, with a finite per-ticket KL bound uniform over
that compact answer interval. The inclusion $A_\theta\subseteq
[-\varepsilon,\varepsilon]$ follows because every fiber contains zero.

\subsection{Characteristic cost and an attaining construction}
\begin{theorem}[Complete characteristic-cost characterization]
\label{thm:fullgame}
For every fixed interior strict point, positive costs and
$0\le\rho<\varepsilon<1$ with $\rho<1/2$,
\begin{equation}\label{eq:fulloptimal}
 \lim_{\delta\downarrow0}
 \frac{B^*(\theta;\varepsilon,\delta)}{\log(1/\delta)}
                         =\frac1{D^*(\theta)}.
\end{equation}
There is one globally honest, independent-pilot family attaining this limit
at all such points, without knowing the true signal vector, the cancellation
stratum, or a unique optimal answer. Its construction is mathematical: it
uses global infima and a growing finite design catalog. The two-ticket
calculation below makes the fixed-instance game certifiably computable;
no polynomial-time general-dimensional optimizer is asserted.
\end{theorem}

\begin{proof}
\textbf{Lower bound for arbitrary real outputs.}
Fix $\theta$ and $\zeta>0$. The two base-fiber endpoints imply
\[
 \Pp_\theta(\tau<\infty,\widehat F\notin A_\theta)\le2\delta.
\]
For each $(a,w)\in A_\theta\times\Delta_r$, select an interior
$x\in\mathcal B(a)$ with
$\sum_j(w_j/c_j)\kl(\theta_j,x_j)<D^*+\zeta$.
Select an actual gate endpoint at $x$ that makes $a$ erroneous. Its error
margin is strict, so it also excludes all outputs in a neighborhood of $a$.
Continuity in $w$ provides a neighborhood of $w$ on which its KL cost is at
most $D^*+2\zeta$. Compactness of $A_\theta\times\Delta_r$ gives a finite
cover and finitely many such alternatives, indexed by $h\le M$.

At a finite return, the realized cost fractions are
$W_j=c_jN_j/C_\tau$. For the covering alternative associated with
$(\widehat F,W)$, the log-likelihood drift is at most
$(D^*+2\zeta)C_\tau$. Under the base law, the centered likelihood ratio
against each fixed interior alternative is a martingale with conditional
variance at most $V$ per cost unit, for a finite $V$ common to the cover.
Freeze the experiment before exceeding a cost cap $T$. Doob's inequality
and change of measure on the excluded output neighborhood give
\[
 \Pp_\theta(C_\tau\le T)
 \le2\delta+M\delta\exp\{(D^*+2\zeta)T+b\}+MV T/b^2.
\]
A zero-query return is covered by assigning it an arbitrary $W$.
Set $\ell=\log(1/\delta)$,
$T=(1-\xi)\ell/(D^*+2\zeta)$ and $b=\xi\ell/2$.
The right side tends to zero. Hence
$\E C_\tau\ge T\Pp(C_\tau>T)$ yields the lower bound, after
$\xi,\zeta\downarrow0$. This finite-cover argument is uniform over algorithms;
it also handles the infimum defining $B^*$ and does not require termination
at alternatives.

\textbf{A small continuity fact.}
For fixed positive $w$, the infimum defining $d$ is continuous in the
interior base parameter. A fixed bad interior point bounds it above, while
positive weights force any relevant minimizer away from the faces of the
parameter cube. Uniform continuity on the remaining compact set proves the
assertion. The same value is continuous in $a$ on the accuracy interval:
the two bad sets are product sublevel and superlevel sets, their interior
threshold level sets are approachable from either side, and thresholds at
$\pm1$ have divergent KL under positive weights. At threshold zero, the
level set is the union of the coordinate hyperplanes $x_j=1/2$ and is
approachable from both product signs. Taking a minimum over the two sides
preserves continuity; they cannot both disappear because $\varepsilon<1$.
A boundary design is approximable from inside the simplex since
\[
 d(\theta,(1-\kappa)w+\kappa\mathbf1/r,a)
                         \ge(1-\kappa)d(\theta,w,a).
\]
At a saturated threshold tending to $1$ or $-1$, every magnitude in
the corresponding product boundary tends to one. With positive weights
and an interior base, that side has divergent KL; the other side is finite
and controls the minimum. At zero, approaching any point on the union of
coordinate hyperplanes by nonzero products requires only a perturbation
tending to zero. These observations justify continuity at the two exceptional
threshold types as well.
Consequently positive designs and a dense answer catalog, including
$\pm\varepsilon$, approximate $D^*$ arbitrarily well.

\textbf{A universal pilot construction.}
Write $L=\max\{1,\ell\}$, $n_0=\lceil L^{2/3}\rceil$,
$h=L^{-1/8}$ and $\kappa=L^{-1/8}$. Take $n_0$ observations per ticket;
clip means into $[1/(2(n_0+1)),1-1/(2(n_0+1))]$. Around the result form a
closed radius-$h$ box $H_L$, intersected with the open cube using half that
clipping distance. On a good pilot event of probability at least
$1-2r\exp(-bL^{5/12})$, the box contains $\theta$ and shrinks to it.
The constant $b>0$ is valid for all sufficiently large $L$ at a fixed interior
base point, by Hoeffding and the fact that clipping is then inactive.

Use finite answer grids of mesh tending to zero in
$[-\varepsilon,\varepsilon]$, always containing the two endpoints, and
finite cost-weight grids becoming dense in the simplex. Mix every grid weight
with $\kappa\mathbf1/r$. Select a catalog pair maximizing the robust score
\[
 S_L(a,w)=\inf_{v\in H_L}d(v,w,a).
\]
If all scores vanish, choose an arbitrary pair. Global exact maximization over
this finite catalog suffices; errors tending to zero would also suffice.
For every fixed $\zeta>0$, the continuity fact and a fixed near-optimal
positive catalog limit show, uniformly on good pilots for all large $L$,
\[
 d(\theta,\bar w,\bar a)\ge S_L(\bar a,\bar w)\ge D^*-\zeta.
\]
No convergence of the selected answer or uniqueness of a maximizer is needed:
the selected pair is frozen before collecting fresh main-stage data.

\textbf{Honest certification.}
Conditionally on the pilot, take cumulative counts
$n_j(k)=\lceil\alpha_j k^2\rceil$, where
$\alpha_j=(\bar w_j/c_j)/\sum_i\bar w_i/c_i$.
Certify $\bar a$ only if
\[
 Z_{\bar a}=\inf_{x\in\mathcal B(\bar a)}
             \sum_j n_j\kl(\widehat\theta_j,x_j) > \beta.
\]
Here $\beta=\ell+\sum_j\log(n_j+1)$ is a valid anytime boundary by
Proposition~\ref{prop:mixture}; its sharper data-dependent version may be used.
Add a fallback: project the rectangle with radii
$\sqrt{\beta/(2n_j)}$ through $[L_\rho,U_\rho]$ and return the midpoint
when its width is at most $2\varepsilon$.
On the simultaneous coverage event the true parameter belongs to the KL
region and the rectangle. Both tests are therefore honest. Conditioning and
integrating preserve the error budget, since all main data are fresh.
For a true signal vector on the boundary of $[0,1]^r$, any strictly
erroneous answer has erroneous interior alternatives approaching that vector,
by continuity of the fiber endpoints. On observations in the support of the
true law, their KL objectives converge to the objective at the true vector.
The infimum over the interior bad set is therefore no larger than that
objective. The same coverage argument proves honesty at endpoints as well;
the efficiency statement continues to require an interior strict instance.

\textbf{Expected cost, including unsuccessful pilots.}
At a fixed strict point there is $u>0$ for which a rectangle contained in a
$2u$-neighborhood has projected width less than $2\varepsilon$.
Every design has sample fractions at least a constant times $\kappa$.
For $k^2\ge K\kappa^{-1}[L+\log(1/\kappa)+\log(L+1)]$, all fallback
radii are at most $u$. Nontermination then requires an empirical deviation
larger than $u$, with probability at most $2r e^{-a\kappa k^2}$.
Summing cost increments $O(k+1)$ times this tail proves the pilot-uniform bound
\begin{equation}\label{eq:newbadpilot}
 \E[C_{\rm main}\mid\text{pilot}]
 \le K'\kappa^{-1}[L+\log(1/\kappa)+\log(L+1)].
\end{equation}
It also proves finite expected termination at every strict point.

For good pilots use the information lower bound $D^*-\zeta$.
Uniformly in alternatives $x$, weights $w$, and nearby $p$, direct expansion
of Bernoulli cross entropy on an interior neighborhood gives
\[
 \sum_j\frac{w_j}{c_j}\kl(p_j,x_j)
 \ge(1-Ku)\sum_j\frac{w_j}{c_j}\kl(\theta_j,x_j)-Ku
 \quad\text{when }\|p-\theta\|_\infty\le u.
\]
To see this, use $p_j/\theta_j,(1-p_j)/(1-\theta_j)\ge1-O(u)$ in the
nonnegative cross-entropy terms; the negative entropy terms are Lipschitz
on that neighborhood. Thus one fixed sufficiently small $u$ makes the GLR
information at least $D^*-2\zeta$ for all selected good-pilot pairs.
Integer rounding contributes only a bounded cost: if $C_k=\sum c_jn_j(k)$,
then $n_j(k)\ge(C_k-\sum c_i)\bar w_j/c_j$.
After $C_k\ge(1+\xi)\ell/(D^*-2\zeta)$ the test has stopped unless the
empirical means deviate by more than $u$. The uniform tail is
$2r e^{-a\kappa C_k}$; at these costs $\kappa L=L^{7/8}\to\infty$.
Its integrated contribution is $o(L)$, and the stage overshoot is
$O(\sqrt L)$. The bad-pilot contribution is $o(L)$ by
\eqref{eq:newbadpilot} and its exponentially small probability.
Pilot cost is $n_0\sum c_j=o(L)$. Let $\xi,\zeta\downarrow0$ to obtain
$\limsup\E C_\tau/\ell\le1/D^*$. Together with the lower bound this
proves the theorem.
\end{proof}

This theorem closes the \emph{variational} finite-slack optimum in the specified
model. It is not a closed-form scalar formula in every dimension. The
continuous-answer principle overlaps directly with existing pure-exploration
theory; the following calculations are what make this model concrete.

\section{Explicit information at cancellation}
\label{sec:explicitnew}
The general game optimizes over answers, allocations and observable
alternatives. At cancellation, symmetry removes the answer search, while a
single increasing scalar equation replaces the allocation game. The inverse
KL magnitude below is the signed-mean magnitude attained at a specified
information level, on the branch above the base magnitude. Equalizing the
coordinate information costs produces a candidate alternative; convexity
in log magnitudes then certifies both that alternative and its allocation.

Write $m_j=|2\theta_j-1|$ and
$k(m,v)=\kl((1+m)/2,(1+v)/2)$. For $D>0$ let $v_j(D)\in(m_j,1)$
solve $k(m_j,v_j(D))=c_jD$. Set $g=(\varepsilon-\rho)/(1-\rho)$.

\begin{theorem}[Inverse-KL product formula]\label{thm:generalgame}
If $|P(\theta)|<g<1$, the information game for the answer zero is the
unique positive solution $D$ of
\begin{equation}\label{eq:newroot}
 \prod_j v_j(D)=g.
\end{equation}
An optimal allocation of cost is
\begin{equation}\label{eq:newweights}
 w_j^*=\frac{c_j/s_j}{\sum_i c_i/s_i},\qquad
 s_j=\frac{v_j(D)(v_j(D)-m_j)}{1-v_j(D)^2}.
\end{equation}
If $P(\theta)=0$, zero is an optimal answer among all real outputs, so
$D^*(\theta)=D$. At the fully fair vector this becomes
\[
 \prod_j(1-e^{-2c_jD})=g^2,\qquad
 w_j^*=\frac{c_j/(e^{2c_jD}-1)}{\sum_i c_i/(e^{2c_iD}-1)}.
\]
\end{theorem}
\begin{proof}
Zero fails exactly when $|P(x)|>g$. For prescribed magnitudes, matching each
alternative sign to the base sign minimizes its KL. Raising an alternative
magnitude below $m_j$ to $m_j$ reduces KL and preserves the constraint.
Thus it is enough to minimize over $v_j\ge m_j$ with $\prod v_j\ge g$.
The left side of \eqref{eq:newroot} is continuous and strictly increases
from $\prod m_j$ to one, proving existence and uniqueness.

Set $z_j=\log v_j$. On the relevant domain,
\[
 \frac{d}{dz}k(m,e^z)=\frac{v(v-m)}{1-v^2},\quad
 \frac{d^2}{dz^2}k(m,e^z)=\frac{v(2v-m-mv^2)}{(1-v^2)^2}>0.
\]
For the displayed weights the derivatives of
$\sum(w_j^*/c_j)k(m_j,e^{z_j})$ at the root vector are equal and positive.
The convex tangent inequality on $\sum z_j\ge\log g$ proves that this
vector minimizes the objective and has value $D$. Conversely the same vector
has value $D$ for every $w$, because each coordinate KL is $c_jD$.
This supplies a direct saddle certificate without a minimax interchange.

If some coordinate is fair, reversing its alternative sign preserves every
coordinate KL and reverses the product. For any nonzero answer $a$, choose
the opposite product sign. The extreme target becomes erroneous at a magnitude
no larger than the threshold $g$ for answer zero. Consequently
$d(\theta,w,a)\le d(\theta,w,0)$ for every $w$. The full game is therefore
solved by zero. Finally $k(0,v)=-\tfrac12\log(1-v^2)$ gives the fully fair
formulas.
\end{proof}

\paragraph{Cost weights versus query fractions.}
The weights in \eqref{eq:newweights} are fractions of expenditure, not of
observations. Query fractions are
$\alpha_j=(w_j/c_j)/\sum_i(w_i/c_i)$.
At $\theta=(.5,.8)$, $c=(1,4)$ and $g=.01$, numerical solution gives
$D\simeq1.27943\,10^{-4}$ and $w\simeq(.96190,.03810)$,
but $\alpha\simeq(.99019,.00981)$.
Thus about 3.81\% of expenditure buys only 0.98\% of queries to ticket two.
These rounded values illustrate the distinction; they are not a new global
certificate or an estimate learned without pilot data.

\paragraph{Connection with classical variance allocation.}
This connection describes an overlap with established product estimation,
not an additional novelty claim. Fix an interior vector with all $m_j>0$,
put $p_0=\prod_jm_j$, and let the zero-answer threshold be
$g=p_0+h$, $h\downarrow0$. Define
\[
 A_j=(1-m_j^2)\prod_{i\ne j}m_i^2,\qquad
 V_c=\left(\sum_j\sqrt{c_jA_j}\right)^2.
\]
Then the zero-answer information $D_0$ and its cost weights satisfy
\[
 D_0\sim\frac{h^2}{2V_c},\qquad
 w_j^*\longrightarrow\frac{\sqrt{c_jA_j}}{\sum_i\sqrt{c_iA_i}}.
\]
Indeed, $v_j=m_j+\sqrt{2c_j(1-m_j^2)D_0}+O(D_0)$ gives
$h=\sqrt{2D_0}\sum_j\sqrt{c_jA_j}+O(D_0)$; substitution in
\eqref{eq:newweights} gives the weights.
For independent signed observations $Y_j=2X_j-1$ and deterministic
counts $n_j$, the exact product-estimator variance is
\[
 \operatorname{Var}\left(\prod_j\overline Y_j\right)
 =\prod_j\left(m_j^2+\frac{1-m_j^2}{n_j}\right)-p_0^2.
\]
At cost $C$ with $n_j=Cw_j/c_j+O(1)$ and fixed positive $w$, its
leading term is $C^{-1}\sum_jc_jA_j/w_j$. Cauchy--Schwarz minimizes
this at $V_c/C$ with the same cost weights; query fractions are proportional
to $\sqrt{A_j/c_j}$. For uncentered Bernoulli means and unit costs this
reduces to the coefficient-of-variation allocation in
\citet{benkamra2015}, equation (11). This deterministic-count identity
alone gives neither anytime coverage nor an adaptive expected stopping cost.
The expansion above fixes nonzero magnitudes and is not uniform as a factor
vanishes; $D_0$ is not asserted to be the full answer-optimized game away
from cancellation. The exact formula, rather than this regular expansion,
underlies the following cancellation limits.

\begin{corollary}[Order and price of cancellation]\label{cor:newstrata}
At a fixed interior vector with $m\ge1$ fair coordinates $Z$, let
$B=\prod_{j\notin Z}|2\theta_j-1|$ and
$G_Z=(\prod_{j\in Z}c_j)^{1/m}$. As $g\downarrow0$,
\[
 D^*(\theta)\sim\frac{g^{2/m}}{2G_ZB^{2/m}},\qquad
 w_j^*\to1/m\ (j\in Z),\qquad w_j^*=O(g^{1/m})\ (j\notin Z).
\]
Thus the high-confidence cost constant grows as $g^{-2/m}$, with the
confidence limit taken first. These are local, not global minimax, rates.
\end{corollary}
\begin{proof}
For a fair coordinate $v_j\sim\sqrt{2c_jD}$. For a fixed biased coordinate,
\[v_j=m_j+\sqrt{2c_j(1-m_j^2)D}+o(\sqrt D).\]
Substitution into \eqref{eq:newroot} gives
$g\sim B(2D)^{m/2}(\prod_{j\in Z}c_j)^{1/2}$.
In \eqref{eq:newweights} the unnormalized weights are $1/(2D)(1+o(1))$
on $Z$ and $O(D^{-1/2})$ elsewhere. Normalization gives the result.
\end{proof}

The same root equation also describes moving cancellation strata. Suppose at
least one coordinate remains fair and the others satisfy
$m_j(g)\sim a_jg^{\beta_j}$, with $\beta_j\ge0$, $a_j>0$, and
$a_j\in(0,1)$ when $\beta_j=0$.
Assign $\beta_j=\infty$ to fair coordinates. The unique solution of
$\sum_j\min(\beta_j,\gamma)=1$ determines the exponent:
$D\sim d g^{2\gamma}$, where
\[
 \left(\prod_{\beta_j<\gamma}a_j\right)
 \left(\prod_{\beta_j=\gamma}(a_j+\sqrt{2c_jd})\right)
 \left(\prod_{\beta_j>\gamma}\sqrt{2c_jd}\right)=1.
\]
This follows by substituting $D=tg^{2\gamma}$ into the inverse-KL expansions
and bracketing the unique increasing product root. The factor with
$\beta_j>\gamma$ is nonempty, so its coefficient equation increases from zero
to infinity. Full moving-stratum weights are derived in Appendix~\ref{app:uniform}. This algebraic statement does not interchange
confidence and precision limits.

\section{Certified computation for two tickets}
\label{sec:quarticnew}
For $a\in A_\theta$, define $b_-(a),b_+(a)$ as the product thresholds of its
lower and upper bad sets. With $x=a-\varepsilon\le0$ and
$y=a+\varepsilon\ge0$, exact sharing gives
$b_- =\max(-1,x)$ and $b_+=\min(1,y)$.
For $\rho>0$, the thresholds are
\[
 b_- =\begin{cases}
 -1,&x\le-1,\\
 (x+\rho)/(1-\rho),&-1<x<-\rho,\\
 1+x/\rho,&-\rho\le x\le0,
 \end{cases}
 \qquad
 b_+ =\begin{cases}
 y/\rho-1,&0\le y<\rho,\\
 (y-\rho)/(1-\rho),&\rho\le y<1,\\
 1,&y\ge1.
 \end{cases}
\]
A threshold $b_-=-1$ or $b_+=1$ has no bad alternatives on that side
and is omitted. Otherwise the bad sets are $P<b_-$ and $P>b_+$.

\begin{theorem}[Quartic projection and a global certificate]\label{thm:quartic}
Let $r=2$, $\theta\in(0,1)^2$ and $0<w<1$. Put
$m_j=2\theta_j-1$, $q_1=w/c_1$, $q_2=(1-w)/c_2$.
For $0<|b|<1$, the information projection onto $P=b$ is the minimum of
\[
 q_1 k(m_1,t)+q_2 k(m_2,b/t)
\]
over all real roots, in $(-1,-|b|)\cup(|b|,1)$, of
\begin{equation}\label{eq:quartic}
 q_1(t-m_1)t(t^2-b^2)-q_2b(b-m_2t)(1-t^2)=0.
\end{equation}
Here $k$ is evaluated for signed arguments in $(-1,1)$.
At $b=0$ the value is
\[\min\{q_1\kl(\theta_1,1/2),q_2\kl(\theta_2,1/2)\}.\]
When $b_-<P(\theta)<b_+$, $d(\theta,w,a)$ is the minimum of the
non-omitted boundary projections; otherwise it is zero.

For rational inputs, root isolation by Sturm sequences and rational logarithm
bounds produce rigorous lower and upper bounds. A finite partition of
$A_\theta$ and feasible alternative cuts gives a global upper bound on $D^*$.
Together with a certified candidate this yields a checkable global optimality
gap, without trusting a floating-point optimizer.
\end{theorem}
\begin{proof}
If the base product lies strictly inside its good band, the minimum over a
bad halfspace occurs on its threshold boundary. To see this, move any strict
bad alternative toward the base along a line segment in parameter space.
Continuity of the product gives a first crossing, while convexity of KL in
its second argument decreases the objective. On $t_1t_2=b\ne0$ write
$t_1=t,t_2=b/t$. The objective diverges at both ends of each admissible interval,
since both weights are positive and both base means are interior.
Every minimum is therefore stationary. Its derivative is
\[
 q_1\frac{t-m_1}{1-t^2}
 -q_2\frac{b(b-m_2t)}{t(t^2-b^2)},
\]
which gives \eqref{eq:quartic}; the denominators do not vanish on the two
intervals. At product zero at least one coordinate equals zero in signed
coordinates; keep the other at its base value. If the base belongs to the
closure of a bad set, its information distance is zero.

For rational inputs the polynomial has rational coefficients. Sturm counts
isolate every distinct real root, including multiple roots. On a root bracket,
$\kl(p,z)$ is bounded using its minimum at the projection of $p$ onto the
bracket and its maximum at an endpoint. Apply the same construction to $b/t$.
To enclose logarithms, reduce a positive rational number to $2^k y$ with
$1\le y<2$, put $z=(y-1)/(y+1)$, and use
\[
 2\sum_{i=0}^{M-1}\frac{z^{2i+1}}{2i+1}
 \le\log y\le
 2\sum_{i=0}^{M-1}\frac{z^{2i+1}}{2i+1}
 +\frac{2z^{2M+1}}{(2M+1)(1-z^2)}.
\]
Use the same formula for $\log2$ and reverse bounds when multiplying by
negative coefficients. Taking the minimum over all isolated roots bounds the
global projection. The implementation rounds final rational enclosures
outward to dyadic rationals to limit integer growth.

For the outer game partition $A_\theta$ into intervals $[a_l,a_u]$.
A lower-side cut with product at most $b_-(a_l)$, or an upper-side cut with
product at least $b_+(a_u)$, is in the closure of the bad set for \emph{every}
answer in that interval. Omitted sides are never used. Choose interior
rational cuts of these types, with per-ticket information bounded above by
$v_{h1},v_{h2}$. Then, throughout this answer interval,
\[
 d(\theta,w,a)\le \min_h\{w v_{h1}+(1-w)v_{h2}\}.
\]
The maximum of this concave piecewise-linear envelope over $w\in[0,1]$
is attained at an endpoint or an intersection of two lines, all computable
rationally. Maximize these bounds over all answer intervals for a global
upper certificate. A certified positive-weight candidate gives a global lower
certificate. Cuts can be generated numerically, but their coordinates, product
constraints, KL upper bounds and coverage of the answer interval are all
rechecked exactly. No unverified minimax equality is used.
\end{proof}

The supplied certificate includes the candidate, answer partition and every
cut. Its verifier recomputes the bounds. Refining the partition and adding cuts
can improve the bracket; the implementation has a work cap and reports the
actual gap even if its requested tolerance is not reached. This is not a
uniform bit-complexity bound. The empirical stopping validator conservatively
abstains at empirical Bernoulli endpoints, which by itself does not affect ideal pointwise interior efficiency.
The executable numerical policy has the precision qualifications stated below.

\section{A smaller anytime boundary}
\label{sec:mixturenew}
An ordinary fixed-sample confidence interval need not retain its coverage
when inspected repeatedly and used to decide when to stop. Here the
confidence event holds simultaneously over the observations collected by
the auditor. On that one event, any answer whose bad set has been excluded
is valid, even if the answer was selected after seeing the data.
The resulting likelihood-based exclusion statistic is a generalized
likelihood ratio (GLR). It controls the decision to return; the pilot's
allocation estimate controls how quickly enough evidence is collected.

The earlier prototype used the method-of-types boundary
$\ell+(r+2)\log(n+1)+\log(\pi^2/6)$. Its finite-confidence overhead was
large. The following replacement is an application of classical mixture
martingales, not a novel concentration inequality.

\begin{proposition}[Beta-binomial simultaneous region]\label{prop:mixture}
Let $N_j,S_j$ be the counts and successes under any predictable ticket
selection. With $\widehat\theta_j=S_j/N_j$ and zero contribution when $N_j=0$,
define
\[
 \mathcal R(N,S)=\sum_{j:N_j>0}\left[
 \log\{(N_j+1)\binom{N_j}{S_j}\}
 +S_j\log\widehat\theta_j+(N_j-S_j)\log(1-\widehat\theta_j)\right].
\]
Use $0\log0=0$. With probability at least $1-\delta$, simultaneously at all
times,
\begin{equation}\label{eq:mixregion}
 \sum_jN_j\kl(\widehat\theta_j,\theta_j)
             \le\ell+\mathcal R(N,S).
\end{equation}
Moreover $0\le\mathcal R\le\sum_j\log(N_j+1)$.
The resulting region permits adaptive designs and data-dependent choice of
answer without an additional union bound over answers or checking times.
\end{proposition}
\begin{proof}
For a fixed alternative vector $x$, the sequential likelihood ratio against
$\theta$ is a nonnegative martingale even with predictable ticket selection.
Mix it over independent uniform priors $x_j\in[0,1]$. Tonelli's theorem gives
a nonnegative martingale starting at one:
\[
 M_\theta=\prod_j
 \frac{\mathrm B(S_j+1,N_j-S_j+1)}
 {\theta_j^{S_j}(1-\theta_j)^{N_j-S_j}}.
\]
The beta integral is $[(N_j+1)\binom{N_j}{S_j}]^{-1}$.
Taking logs gives
$\log M_\theta=\sum_jN_j\kl(\widehat\theta_j,\theta_j)-\mathcal R$.
Ville's inequality proves \eqref{eq:mixregion}. The binomial mass at its
maximum-likelihood parameter lies between $1/(N_j+1)$ and one. For the lower
bound on that mass, $S_j$ is a mode of $\operatorname{Bin}(N_j,S_j/N_j)$
and there are $N_j+1$ possible counts. This proves both regret bounds.
Any answer whose entire bad set lies outside this region is simultaneously
valid; data dependence of its selection cannot remove the true point from
the region. Endpoint true parameters follow by the same mixture
supermartingale argument on their support, or by continuity of bounded
transcripts.
\end{proof}

At fixed interior empirical proportions, Stirling's formula gives
$\mathcal R=\frac12\sum_j\log N_j+O(1)$, rather than the previous
$(r+2)\log n$ overhead. This is an improvement over our previous boundary,
not a comparison with the strongest existing mixture boundaries, some of which
use iterated logarithms. Exact rational bounds are supplied for moderate
integer counts; floating gamma evaluations in simulations are labelled as
numerical diagnostics. The looser bound $\ell+\sum_j\log(N_j+1)$ remains a
simple certified option for very large counts.

\section{Joint limits: a sufficient regime and a counterexample}
\label{sec:jointnew}
This section concerns the fixed-pilot zero-certificate auditor, rather than the
catalog construction of Theorem~\ref{thm:fullgame}. Its pilot takes
$n_0=\lceil L^{2/3}\rceil$ observations per ticket, where
$L=\log(1/\delta)$. If the pilot product magnitude is below $g$, it uses
\eqref{eq:newweights} at the pilot means; otherwise it uses uniform cost
weights. It mixes the selected weights with $L^{-1/8}\mathbf1/r$ and uses
fresh main observations, a zero GLR, and the rectangular fallback. The older sufficient joint-limit regime and its technical stability proof
are retained in Appendix~\ref{app:oldjoint}. Here we isolate the failure that
motivates the new pilot.

\begin{proposition}[A factor-two failure outside the resolved-pilot regime]
\label{prop:pilotfailure}
Take $r=2$, unit costs, $\theta=(1/2,4/5)$, $g_L=1/L$, and fixed
$0\le\rho<1/2$. For the preceding fixed-pilot auditor, with its original type
boundary or the beta-binomial boundary (and either the quadratic or
vanishing-ratio geometric checking schedule),
\[
 \liminf_{L\to\infty}\frac{D_L\E_\theta C_\tau}{L}\ge2.
\]
Thus its pointwise fixed-precision optimality cannot be promoted to an
unrestricted simultaneous confidence--precision claim.
\end{proposition}
\begin{proof}
Write $B=3/5$. A pilot fair-coordinate magnitude is typically of order
$n_0^{-1/2}=L^{-1/3}$, while the biased-coordinate magnitude converges to $B$.
More explicitly, outside an exponentially unlikely biased-coordinate event,
selecting the inverse-KL branch requires the fair pilot mean to lie within
$O(1/L)$ of $1/2$. The interval contains at most a bounded number of binomial
counts, each having probability $O(n_0^{-1/2})$. Consequently the uniform-cost
fallback is selected with probability $1-o(1)$. Its regularization remains
exactly uniform.

Condition on that design and cap the main cost at a fixed multiple of $L^3$.
A method-of-types union bound over the polynomially many possible stage totals
gives, with probability $1-o(1)$,
$\sum_j n_j\kl(\widehat\theta_j,\theta_j)\le K\log L$ at every checked
stage below this cap. For the zero GLR use an alternative with second mean
$4/5$ and first signed magnitude $g_L/B$, matching the empirical sign of the
first coordinate. The identity
\[
 k(m,v)=k(m,0)+k(0,v)-m\operatorname{arctanh}(v)
                  \le k(m,0)+k(0,v),\quad m,v\ge0,
\]
implies
$Z_g\le K\log L+n_1[-\tfrac12\log(1-(g_L/B)^2)]$.
As the boundary is at least $L$, this test cannot stop before
$n_1\ge(1-o(1))2B^2L/g_L^2$.

The fallback cannot stop earlier either. Its first radius is at least
$\sqrt{L/(2n_1)}$, whereas the empirical error is only
$O(\sqrt{\log L/n_1})$. The rectangle contains a symmetric interval around
the fair mean with radius $(1-o(1))\sqrt{L/(2n_1)}$, clipped to $[0,1]$,
and contains the true second coordinate. Its product range therefore contains
$\pm B\min\{1,(1-o(1))\sqrt{2L/n_1}\}$.
The projected robust width can be at most $2\varepsilon_L$ only if
$n_1\ge(1-o(1))2B^2L/g_L^2$. Under uniform sampling, total main cost is
$2n_1+O(1)$, so with probability tending to one it is at least
$(1-o(1))4B^2L^3$. Meanwhile $D_L\sim1/(2B^2L^2)$ by
Corollary~\ref{cor:newstrata}. The expected-cost lower limit follows from
nonnegativity. The argument uses only a boundary at least $L$ and polynomially many
checks below the cap, so applies to both boundaries and schedules stated.
This proves a failure of efficiency, not of honesty.
\end{proof}

\section{Learning the design from charged observations}
\label{sec:learned}
The preceding characteristic cost is useful only if unknown signals can be
handled. This section gives a simpler implementable specialization of the
frozen-pilot construction, and a replacement for the pilot that failed in
Proposition~\ref{prop:pilotfailure}. Neither policy receives the true means.

\begin{lemma}[Balancing the two bad sets]\label{lem:balance}
Fix an interior two-ticket vector, positive weights, and a strictly feasible
precision. Let $d_-(a)$ and $d_+(a)$ be the lower and upper bad-set information
distances, assigning $+\infty$ to an absent bad set. On $A_\theta$,
$d_-$ is nonincreasing and $d_+$ is nondecreasing. A maximizer of
$\min\{d_-(a),d_+(a)\}$ can be found by bisection of their ordering, also
checking the two endpoints of $A_\theta$. Only the remaining weight search
requires global maximization.
\end{lemma}
\begin{proof}
As $a$ increases, $\{L(P)<a-\varepsilon\}$ expands and
$\{U(P)>a+\varepsilon\}$ contracts. Taking infima proves monotonicity.
Before a crossing, the smaller nondecreasing function controls the minimum;
after it, the smaller nonincreasing function controls it. Continuity at finite
values follows from the compact-coercivity argument in
Theorem~\ref{thm:fullgame}. An absent side can occur at a saturated endpoint,
which is explicitly checked. A flat crossing causes no difficulty. The
projection on each side is computed by Theorem~\ref{thm:quartic}.
This argument does not assert concavity of the optimized value in the weight.
\end{proof}

\paragraph{A practical frozen-pilot policy.}
Take $n_0=\max\{32,\lceil L^{2/3}\rceil\}$ observations from each ticket
and discard them for main inference. Clip their means as in
Theorem~\ref{thm:fullgame}. On a finite weight catalog, mix weights with
$\kappa\mathbf1/r$, where $\kappa=\min\{.2,L^{-1/8}\}$, and choose the
largest empirical information score. In the two-ticket implementation the
catalog has $\max\{13,\lceil L^{1/8}\rceil\}$ equally spaced weights;
Lemma~\ref{lem:balance} supplies each answer. Local numerical refinement is
allowed only as an improvement over the catalog winner.
Theoretical score and answer errors must tend to zero as $L\to\infty$.
If all scores are zero, use uniform cost weights. Freeze the selected design,
draw fresh data, check the frozen answer and a newly balanced answer, and use
the projected rectangle as a fallback. Main sample totals grow by a factor
$1+\min\{.15,L^{-1/4}\}$, with integer rounding.

\begin{proposition}[Pointwise optimality without the robust-box optimization]
\label{prop:plugin}
In exact arithmetic, or with vanishing score errors, the preceding catalog
policy is globally honest and attains \eqref{eq:fulloptimal} at every fixed
strict interior instance. The same assertion holds with any increasing dense
catalog in fixed dimension. Selecting additional answers from the main data
does not require a new error budget.
\end{proposition}
\begin{proof}
On a compact neighborhood of a fixed $\theta$, expansion of Bernoulli cross
entropy gives, uniformly over $w,a$ and all alternatives,
\[
 d(\theta,w,a)\ge(1-Ku)d(p,w,a)-Ku
 \quad\text{if }\|p-\theta\|_\infty\le u.
\]
The reverse inequality follows on interchanging $p$ and $\theta$.
Take one fixed positive near-optimal design and answer, approximated by the
catalog. Its empirical score tends to its true score. Hence the empirical
catalog winner has true score at least $D^*(\theta)-o(1)$, even if its answer
does not converge. Choosing $u=L^{-1/8}$ bounds bad-pilot probability by
$2r\exp(-bL^{5/12})$. The uniform bad-pilot fallback bound
\eqref{eq:newbadpilot} remains valid. For good pilots the selected frozen
answer alone gives the expected-cost argument of Theorem~\ref{thm:fullgame}.
Geometric stage overshoot is $o(L)$ and the pilot costs $o(L)$.
Both the frozen and newly balanced answers are checked against the same
simultaneous mixture region, so extra candidate selection cannot compromise
honesty or increase the stopping cost. The general lower bound concludes.
\end{proof}

\begin{corollary}[Reusing the fixed pilot]\label{cor:reuse}
The fixed-pilot policy remains globally honest and pointwise optimal when
pilot observations are retained in the mixture region, charged once, and main
counts are $n_j(N)=\max\{n_0,\lceil\alpha_jN\rceil\}$.
\end{corollary}
\begin{proof}
Start the mixture martingale before the pilot. Subsequent ticket choices are
predictable, so its coverage includes the pilot and main observations even
though the design is estimated from the same data. No conditional-confidence
argument is used for validity.
Total cost at virtual total $N$ is at most
$n_0\sum c_j+N\sum c_j\alpha_j+O(1)$; the extra term is $o(L)$.
The information statistic has counts at least $\alpha_jN$.
For each ticket, the exponential Hoeffding martingale gives
\[
 \Pp\{\exists n\ge m:|\widehat\theta_{j,n}-\theta_j|>u\}
                          \le2e^{-2mu^2}.
\]
Thus the good-pilot cost argument holds without conditioning on those
observations, by bounding the intersection with this simultaneous prefix
event. On a bad pilot, once $N\ge K n_0/(\kappa u)$, its possible contribution
to any pooled mean is at most $u$, and fresh main observations control the
remaining deviation. The conditional fallback cost is therefore at most
$K'\kappa^{-1}(L+n_0)$ at a fixed strict point. Multiplication by the bad-pilot
probability makes its contribution $o(L)$. The same lower bound concludes.
\end{proof}

The implemented floating optimizer is a proposal engine, not a claim of
machine-certified asymptotic efficiency. Every reported experimental return
is authorized again by rational inequalities. Their validity does not depend
on the design proposal finding an optimal weight. The authoritative returned
estimate is stored as a rational number, with a floating display value. The theorem
concerns the uncapped policy; experimental budget caps are reported explicitly.
The revised return guards increase KL rounding precision, logarithm depth,
root-isolation depth and rectangle precision with the observation count.
This removes the fixed absolute grids identified in review; it does not prove
optimality of the floating design proposals. Supported query caps are integers
from zero to $2^{52}$, checked before every charged batch, including pilots.
When an exact risk is supplied it controls both design and inference.
Historical experiment records retain their original numerical implementation.

\paragraph{Fast rigorous boundary evaluation.}
Direct exact binomial coefficients become cumbersome at large counts.
We use the classical Robbins enclosure \citep{robbins1955}
\[
 (n+\tfrac12)\log n-n+\tfrac12\log(2\pi)+\frac1{12n+1}
 <\log(n!)<
 (n+\tfrac12)\log n-n+\tfrac12\log(2\pi)+\frac1{12n},
\]
with $0!=1!=1$ handled exactly. Alternating-series bounds and Machin's
identity enclose $\pi$, and the rational logarithm bounds used earlier enclose
all other terms. Substitution into the beta-binomial regret, always in the
conservative direction, gives a rational upper boundary. Quartic root brackets
give a rational lower GLR. A return is allowed only if the latter strictly
exceeds the former; endpoint empirical means use the rectangular fallback.
These numerical ingredients are classical, not separate novelty claims.

\section{An information-scaled pilot repairs the joint-limit failure}
\label{sec:resolvednew}
For zero certification put $g=(\varepsilon-\rho)/(1-\rho)$ and $L=\log(1/\delta)$.
At pilot round $k\ge0$, sample every ticket to $n_k=2^k$ observations and
form clipped means $\widetilde\theta_k$. If
$|P(\widetilde\theta_k)|<g/2$, compute the inverse-KL product root
$\widehat D_k$ and its cost weights. Leave the pilot at the first round with
\begin{equation}\label{eq:pilottrigger}
 n_k\widehat D_k\ge\Lambda_{L,k},\qquad
 \Lambda_{L,k}=\sqrt L+\log^2(k+2).
\end{equation}
At every pilot round also check the projected confidence rectangle, with
risk $\delta/2$, and return if it certifies a target. This fallback matters:
zero need not be acceptable at every strict fiber. After leaving the pilot,
discard its observations for inference, freeze the computed weights mixed
with $\kappa\mathbf1/r$, and use fresh observations to check zero and the
rectangle at risk $\delta/2$. Use the preceding vanishing-ratio geometric
schedule. Thus the algorithm remains honest and has finite expected cost at
every fixed strict instance, including away from cancellation.

\begin{theorem}[Joint efficiency with an information-scaled pilot]
\label{thm:adaptivejoint}
Fix $r,c,\rho$ and $a_0\in(0,1/2)$. Let
$\theta_L\in[a_0,1-a_0]^r$, $P(\theta_L)=0$, and $g_L\downarrow0$.
Let $D_L$ be the inverse-KL product root at $(\theta_L,g_L)$.
If
\begin{equation}\label{eq:subexponential}
                 \log(1/D_L)=o(L),
\end{equation}
then the information-scaled pilot auditor, at
$\varepsilon_L=\rho+(1-\rho)g_L$, satisfies
\[
                 D_L\E_{\theta_L} C_\tau/L\longrightarrow1.
\]
Its expected pilot cost is $O(\sqrt L/D_L)$ along these arrays. At every fixed
cancellation stratum, all polynomial precisions $g_L=L^{-b}$, $b>0$, satisfy
the condition. In particular this repairs the $g_L=1/L$ counterexample.
The assertion still does not cover arbitrary joint limits.
\end{theorem}
\begin{proof}
For analysis, realize the potential pilot and main observations on independent
infinite streams. An event involving unused pilot observations then places no
restriction on main data. Constants are uniform on the compact parameter cube. The inverse-root
and weight expansions used in Theorem~\ref{thm:jointnew} remain uniform.
Appendix~\ref{app:uniform} supplies uniform root, weight, GLR and fallback
lemmas and the stopping-time tail sums.
We give the pilot argument explicitly, since a fixed-precision continuity
claim would be insufficient here. For any fixed strict instance, uniform
pilot sampling makes the pilot fallback have finite expected cost; after any
trigger the regularized main fallback does too. The two simultaneous coverage
events, each of risk $\delta/2$, prove global honesty.

\textbf{One event for every pilot round.}
Let $h_{L,k}=L^{1/4}+3\log(k+2)+\log(8r)$.
Hoeffding and a summation over $r$ tickets and $k\ge0$ give an event
$\mathcal E_L$ with
\[
 \Pp(\mathcal E_L^c)\le K e^{-L^{1/4}},\qquad
 \|\widehat\theta_k-\theta_L\|_\infty
       \le \sqrt{h_{L,k}/(2n_k)}\quad\text{for every }k.
\]
At any trigger round,
\[
 \frac{h_{L,k}}{\Lambda_{L,k}}\le K L^{-1/4},
\]
uniformly in $k$: put $x=\log(k+2)$ and use
$x/(\sqrt L+x^2)\le1/(2L^{1/4})$.
The pilot root is uniformly bounded above for small $g_L$. Indeed a fixed
KL level raises each magnitude at least to $(1+m_j)/2$, whose product is
at least $2^{-r}>g_L$. Thus a trigger requires $n_k\ge K\sqrt L$.
On $\mathcal E_L$ the pilot then lies in a fixed interior cube; clipping is
inactive, and
\begin{equation}\label{eq:relativepilot}
 \|\widetilde\theta_k-\theta_L\|_\infty
                    \le K L^{-1/8}\sqrt{\widehat D_k}.
\end{equation}

\textbf{Relative root and weight accuracy at a random trigger.}
For two nearby interior magnitude vectors, their inverse-KL magnitudes at a
fixed small level $d$ differ by at most a constant times the base-vector
distance. Indeed implicit differentiation of $k(m,v)=c_jd$ gives
\[
 \frac{\partial v}{\partial m}=
 \frac{(1-v^2)[\operatorname{arctanh}(v)-\operatorname{arctanh}(m)]}{v-m}
 =1+O(\sqrt d),
\]
uniformly on the compact range, with its continuous value when $d=0$. Each root magnitude is at least $K\sqrt d$.
Consequently \eqref{eq:relativepilot} implies that their root products at
$d=\widehat D_k$ have ratio $1+o(1)$, uniformly over trigger rounds.
The product relation and $g_L\to0$ first force $\widehat D_k\to0$.
At the true vector at least one magnitude is zero; its root contributes
$1/2+o(1)$ to the logarithmic elasticity of the product in $d$, while other
contributions are nonnegative. Bracketing that increasing product proves
$\widehat D_k/D_L\to1$. Substitution in the explicit weights proves
$\widehat w_{k,j}/w_j^*\to1$ for all $j$. This also justifies stability when
the identity or number of fair coordinates changes with $L$.

\textbf{Pilot duration, not just its success probability.}
Let $k_*$ be the first index with
$2^{k_*}D_L\ge K\Lambda_{L,k_*}$, for a sufficiently large fixed $K$.
Under \eqref{eq:subexponential},
$k_*=O(\log(1/D_L)+\log L+1)=O(L)$ and
$\Lambda_{L,k_*}=O(\sqrt L)$.
At this round, any mean error at most $b\sqrt{D_L}$, for a sufficiently
small fixed $b>0$, ensures both $|P(\widetilde\theta_k)|<g_L/2$ and the
trigger inequality. The product claim follows by comparing every perturbed
factor with its inverse root and using a true fair coordinate; decreasing
$b$ makes that factor's ratio arbitrarily small. Root stability then gives
$\widehat D_k\ge D_L/2$. The same implications hold at later rounds.
Thus, even without conditioning on $\mathcal E_L$,
\[
 \Pp(\hbox{pilot unfinished at round }k)
                     \le2r\exp(-b'2^kD_L),\qquad k\ge k_*.
\]
Summing the doubling-stage cost increments against this tail yields
$\E C_{\rm pilot}=O(2^{k_*})=O(\sqrt L/D_L)=o(L/D_L)$.
Earlier fallback return can only reduce this cost.

\textbf{Main-stage cost and bad pilots.}
On $\mathcal E_L$, a triggered design has zero-answer information
$(1-o(1))D_L$, after regularization. The relative GLR stability calculation
in Theorem~\ref{thm:jointnew} applies unchanged. At cost proportional to
$L/D_L$, the beta-binomial boundary is
$L+O(\log L+\log(1/D_L))=L[1+o(1)]$; the additional $\log2$ is negligible.
Mean-deviation tails with weights at least $\kappa/r$ have integral
$o(L/D_L)$ after any cost $(1+\xi)L/D_L$.
The relative geometric overshoot tends to zero.

For every possible pilot outcome, the rectangle fallback has conditional
expected main cost at most
\[
 K\kappa^{-1}D_L^{-1}
       [L+\log(1/D_L)+\log(1/\kappa)+1].
\]
This follows from a cube of radius $b\sqrt{D_L}$ contained in the zero-good
band, as established in Theorem~\ref{thm:jointnew}. Multiply this uniform
bound by $\Pp(\mathcal E_L^c)\le K e^{-L^{1/4}}$ to get $o(L/D_L)$.
Fresh main data justify conditioning on the complete stopped pilot.
Together with the pilot estimate this proves the upper limit.

For the lower limit, the two signed product alternatives and finite-cap
likelihood argument of Theorem~\ref{thm:jointnew} need compactness and
$L\to\infty$, not its stronger pilot-resolution condition. Their variance
per cost is $O(D_L)$, so the remainder at cost $L/D_L$ is $O(1/L)$.
Let their product approach $g_L$ from outside and then let the cost slack
vanish. This gives the reverse limit for every globally honest family.
Finally $D_L\asymp g_L^{2/m}$ at a fixed stratum proves the polynomial claim.
\end{proof}

Condition \eqref{eq:subexponential} is a substantial enlargement of
\eqref{eq:resolvedpilot}; it is a sufficient regime, not a sharp frontier.
Pilot selection depending on estimated difficulty is not itself a new
principle. The specific trigger, product-relative stability argument and
joint-limit guarantee are the candidate contribution here.

\section{Precision-aware pilots and localized certification}
\label{sec:localized}
There are two costs to control when precision and confidence change together:
learning a sufficiently accurate allocation and accumulating evidence for a
valid stopping decision. The calibrated pilot below addresses the first;
the checkpoint boundary addresses the second. Their combination yields
Theorem~\ref{thm:loglogjoint}. Its proof must bound \emph{expected} cost,
including rare inaccurate or unusually long pilots. A high-probability
allocation statement alone would not establish that conclusion.

The logarithmic restriction in Theorem~\ref{thm:adaptivejoint} partly pays
for the chosen confidence region. Iterated-logarithm KL boundaries are already
available \citep{kaufmannkoolen2021,howard2021}. Informative mixing priors and
predictive assistance are also established ideas
\citep{cortinovis2025prior,kilian2026bayes}. We do not claim those principles.
Here a simple checkpoint construction and a precision-aware trigger extend
the product theorem to a much wider array regime. A separate pilot-localized
mixture augments the executable boundary without being needed for that theorem.

\subsection{Checkpoint confidence region}
Readers focused on the main joint-efficiency result can proceed from
Proposition~\ref{prop:checkpoint} directly to the calibrated policy in
Section~\ref{subsec:calibrated}. The next two propositions describe optional
certification channels; the theorem does not require the localized mixture.

\begin{proposition}[Conditional checkpoint region]\label{prop:checkpoint}
Condition on an arbitrary stopped pilot. From independent fresh streams take
positive deterministic counts $n_j(s)$ at checkpoints $s=0,1,\ldots$;
the schedules may depend on that pilot but not on main observations.
Fix $0<\lambda<1$ and a main risk $d\in(0,1)$. Then simultaneously over
all checkpoints, with probability at least $1-d$,
\[
 \sum_j n_j(s)\kl(\widehat\theta_j(s),\theta_j)\le
 b_s(d):=\frac{\log(1/d)+\log[(s+1)(s+2)]
                   +r\log[2/(1-\lambda)]}{\lambda}.
\]
The same guarantee holds unconditionally and allows data-dependent answers
and stopping among these checkpoints.
\end{proposition}
\begin{proof}
For one fixed positive count, the two Bernoulli Chernoff tails imply
$\Pp\{Y>x\}\le2e^{-x}$ for $Y=n\kl(\widehat\theta,\theta)$.
Integration of the tail gives
\[
 \E e^{\lambda Y}
 =1+\lambda\int_0^\infty e^{\lambda x}\Pp(Y>x)\,dx
 \le\frac{1+\lambda}{1-\lambda}\le\frac2{1-\lambda}.
\]
At a fixed checkpoint, fresh ticket streams are conditionally independent,
so their exponential moments multiply. Markov's inequality gives error at
most $d/[(s+1)(s+2)]$. These probabilities sum to $d$ by telescoping.
No independence between checkpoints is needed. Integrate over the pilot.
The proof also covers true endpoint means, where the corresponding $Y$ is zero.
The restriction on the count schedule is essential: the fixed-count
exponential-moment argument does not license arbitrary main-data-driven counts.
\end{proof}

\subsection{Optional localized and pooled channels}
\begin{proposition}[Pilot localization with a protected fallback]
\label{prop:localmixture}
Let the stopped pilot have counts $m_j$ and successes $t_j$.
Conditionally use independent $\operatorname{Beta}(t_j+1,m_j-t_j+1)$
mixing priors for fresh main likelihoods. Put
\begin{align*}
 R_{\rm loc}(n,S)&=\sum_j\{n_j h(\widehat p_j)+\log B_{0j}-\log B_{1j}\},\\
 h(p)&=p\log p+(1-p)\log(1-p),\\
 B_{0j}&=\mathrm B(t_j+1,m_j-t_j+1),\\
 B_{1j}&=\mathrm B(t_j+S_j+1,m_j+n_j-t_j-S_j+1),
\end{align*}
with $0\log0=0$ and zero contribution at $n_j=0$.
The region $\sum n_j\kl(\widehat p_j,\theta_j)\le\log(1/d)+R_{\rm loc}$
has anytime conditional coverage $1-d$. On the frozen checkpoint schedule,
take the minimum of this boundary, the uniform-prior boundary of
Proposition~\ref{prop:mixture}, and $b_s$, each at risk $d/3$.
The resulting region has coverage $1-d$. Its boundary is at most the
uniform-prior boundary at risk $d$ plus $\log3$.
\end{proposition}
\begin{proof}
For every fixed alternative vector $q$, the likelihood ratio of fresh main
data against an interior true vector is a nonnegative martingale conditional
on the pilot. At a true endpoint, define the ratio on the support of the true
law; it is a nonnegative supermartingale, since alternative mass outside that
support can be lost. Integrate against the pilot-measurable probability density.
Its initial value is one. Beta integration gives the displayed regret and
Ville gives coverage; predictable main ticket choices are allowed for this
mixture component. The checkpoint component additionally needs the schedule
condition of Proposition~\ref{prop:checkpoint}. A union bound over the three
regions gives their intersection, equivalently the minimum boundary.
The last inequality follows by retaining just the uniform-prior component.
It is a boundary guarantee, not a uniform stopping-cost improvement.
\end{proof}

Pilot observations in this construction choose a density and are charged,
but are not counted again in the main likelihood. Choosing a density after
seeing the very observations to which it is applied is not covered: with one
fair draw, choosing $\operatorname{Beta}(K+1,1)$ after a success and its
reflection after a failure gives expected likelihood ratio
$2(K+1)/(K+2)>1$. The delivered test suite reproduces this counterexample.
No Bayesian correctness assumption on the mixing prior is required.

\begin{proposition}[Pooled rescue without an additional risk budget]
\label{prop:pooledrescue}
Continue the pilot's uniform-beta likelihood-ratio process across all later
observations, without resetting it. At any completed batch, additionally
check zero and a projected rectangle using this process at risk $\delta/2$.
Keep the base policy's conditional main region at total risk $\delta/2$.
This augmented policy is globally honest. Coupled to its base policy on the
same per-ticket streams and randomization, it never uses more observations
or costs more. If the base policy is asymptotically optimal on an array,
the augmented policy remains so.
\end{proposition}
\begin{proof}
Proposition~\ref{prop:mixture} allows arbitrary predictable ticket choices,
including those selected after the pilot. Its single coverage event at risk
$\delta/2$ now covers both the original pilot fallback and every additional
pooled check. The conditional main coverage event has risk at most $\delta/2$
after integration over the potential pilot. Their union proves honesty;
there is no independence assumption between the two events.
Until an extra return, the observation requests and decisions of the base
algorithm are unchanged. Thus either an extra certificate stops earlier,
or both procedures have exactly the same cost. This proves pathwise cost
dominance, including budget exhaustion. Combine the inherited expected-cost
upper bound with the common lower bound for the last assertion.
\end{proof}

The pooled and localized channels have different roles. The first retains a
single likelihood process from time zero. The second conditions on the pilot
and starts a fresh likelihood. Applying a pilot-chosen density retroactively
to pooled observations would be invalid. The rescue guarantee concerns
observation cost; checking an extra region can increase computation time.

\subsection{Precision-aware policy and joint efficiency}
\label{subsec:calibrated}
\paragraph{A calibrated pilot and checkpoint rule.}
Set $L=\log(1/\delta)$ and, using the requested $g$, define
\begin{align*}
 H&=\sqrt L+\log\{e+\log(1/g)\},\qquad
 q=\max\{1,(L/H)^{1/4}\},\\
 h_k&=H+3\log(k+2)+\log(8r).
\end{align*}
Replace \eqref{eq:pilottrigger} by $n_k\widehat D_k\ge qh_k$, retaining
$|P(\widetilde\theta_k)|<g/2$, doubling pilot counts, and the pilot rectangle
at risk $\delta/2$. With $\widehat w$ the pilot cost weights, freeze
\[
 \kappa=\min\{.2,L^{-1/8}\},\qquad
 w=(1-\kappa)\widehat w+\kappa\mathbf1/r,\qquad
 \alpha_j=\frac{w_j/c_j}{\sum_i w_i/c_i}.
\]
Use independent fresh main streams and the virtual-count calendar
\begin{equation}\label{eq:checkpoint-calendar}
 \eta_L=\min\{.15,L^{-1/4}\},\quad N_0=16,\quad
 N_{s+1}=\lceil(1+\eta_L)N_s\rceil,\quad
 n_j(s)=\lceil\alpha_jN_s\rceil.
\end{equation}
Only the additional observations needed to reach these counts are drawn.
Thus the entire calendar is fixed conditionally on the pilot; rounding is
applied to each virtual total, not recursively to each ticket count. Take
$\lambda=1-1/\max\{2,\lceil L\rceil\}$.
Use the checkpoint region at risk $\delta/2$, or the three-component region
at total main risk $\delta/2$. Both a zero-answer GLR and the projected
rectangle are checked at every main checkpoint.

\begin{theorem}[Joint efficiency under an iterated-logarithm condition]
\label{thm:loglogjoint}
Fix $r\ge1$, $c_j>0$, $0\le\rho<1/2$ and $a_0\in(0,1/2)$.
Let $L=\log(1/\delta)\to\infty$, $\theta_L\in[a_0,1-a_0]^r$
with $P(\theta_L)=0$, and
$\varepsilon_L=\rho+(1-\rho)g_L$. Let $D_L$ be the inverse-KL product root
at $(\theta_L,g_L)$. The calibrated policy with
\eqref{eq:checkpoint-calendar} is
globally honest and has finite mean termination at every fixed strict instance.
Along arrays with $g_L\downarrow0$ and
\begin{equation}\label{eq:loglogregime}
                \log\{e+\log(1/g_L)\}=o(L),
\end{equation}
its uncapped exact-arithmetic version satisfies
\[
             D_L\E C_\tau/L\longrightarrow1,
 \qquad \E C_{\rm pilot}=O(L^{1/4}H_L^{3/4}/D_L)=o(L/D_L).
\]
In particular, this includes $g_L=\exp(-L^a)$ for every fixed $a>0$ and
$g_L=\exp\{-\exp(L^b)\}$ for $0<b<1$.
No claim of a necessary frontier, or of numerical implementation at such
extreme scales, is made.
The same conclusion holds after pooled rescue.
\end{theorem}
\begin{proof}
Write $D=D_L$, $H=H_L$. Uniformly on the interior compact cube, the
product root first tends to zero with $g$: at any fixed positive information
level, every inverse magnitude has a uniform positive lower bound. The
inverse-KL increments are then between constant multiples of $\sqrt D$.
Every inverse magnitude is therefore at least $c\sqrt D$, while one fair
coordinate is at most $C\sqrt D$ and all other magnitudes are at most one.
The product equation implies
\begin{equation}\label{eq:Dgcomparison}
          c'g^2\le D\le C'g^{2/r}.
\end{equation}
Consequently the iterated logarithms of $1/D$ and $1/g$ differ by $O(1)$.
Under \eqref{eq:loglogregime}, $H=o(L)$, $q\to\infty$, and
$qH=L^{1/4}H^{3/4}=o(L)$.

\textbf{Pilot accuracy and duration.}
Hoeffding and summation over tickets and rounds give an event $E$ with
\[
 \Pp(E^c)\le K e^{-H},\qquad
 \|\widehat\theta_k-\theta_L\|_\infty\le\sqrt{h_k/(2n_k)}
 \quad\hbox{for every }k.
\]
At any trigger the error is at most $Kq^{-1/2}\sqrt{\widehat D_k}$.
As in Theorem~\ref{thm:adaptivejoint}, the estimated root is uniformly bounded
above; hence a trigger requires $n_k\ge c qH\to\infty$.
On $E$, the error at any trigger is at most $Kq^{-1/2}$ and therefore tends
to zero. The true means lie in $[a_0,1-a_0]$ and the clipping threshold
tends to zero with $n_k$, so clipping is inactive there eventually.
Lemma~\ref{lem:rootstable} yields, uniformly over
trigger rounds,
\[
 \widehat D_k/D=1+o(1),\qquad \widehat w_{k,j}/w_j^*=1+o(1).
\]
Let $k_*$ be the first round with $2^{k_*}D\ge K qh_{k_*}$ for a
sufficiently large constant $K$. A candidate round of order
$\log(1/D)+\log L+1$ shows that
$\log(k_*+2)=O(H)$, using \eqref{eq:Dgcomparison}; minimality and doubling
then give $2^{k_*}=O(qH/D)$.
At every round $k\ge k_*$, an empirical error at most $b\sqrt D$ implies
the product condition and trigger, by Lemmas~\ref{lem:rootstable} and
\ref{lem:cubetail}. The implication persists because $h_{k+1}/h_k<2$.
Thus $\Pp(\hbox{pilot unfinished at }k)\le2r e^{-a2^kD}$.
Summing doubling increments gives $\E C_{\rm pilot}=O(qH/D)$.
This is an unconditional expectation bound; unusually late pilots are included.

\textbf{The new boundary scale.}
Let $C_s$ denote main cost at checkpoint $s$. Fixed positive costs and the
geometric schedule imply, uniformly over frozen designs,
\[
 \log(s+2)\le K+K\log(L+1)+\log\log(e+C_s).
\]
Indeed total counts are the virtual count plus a bounded rounding term,
cost lies between fixed multiples of that count, and the reciprocal log
growth factor is $O(1+L^{1/4})$. For all $L\ge2$ and all checkpoints,
enlarging a constant $K$ depending only on the fixed dimension and costs gives
the explicit envelope, also at risk $\delta/6$,
\begin{equation}\label{eq:checkpoint-envelope}
 b_s\le B_L(C_s):=L+K[1+\log(L+1)+\log\log(e+C_s)].
\end{equation}
Here $(1-1/\lceil L\rceil)^{-1}L\le L+2$; its remaining bounded
factor is absorbed into $K$. For $C>0$,
$C B_L'(C)\le K$, whereas $B_L(C)\ge L+K$.
Thus $B_L(C)/C$ is decreasing. This is a monotone upper envelope;
the discrete boundary itself need not have a decreasing ratio.
At $C_s\asymp L/D$ it gives
\[
 L+O\{\log(L+1)+\log\log(e+1/D)\}=L+o(L).
\]
This replaces the previous $\log(1/D)$ overhead. The hybrid boundary is
no larger than this checkpoint boundary.

\textbf{Good-pilot expected main cost.}
On $E$, the frozen regularized design has information $(1-o(1))D$.
To verify the needed weight bound, the uniform inverse-root increments give
$cD\le s_j\le C\sqrt D$ in \eqref{eq:newweights}. Consequently
$c_j/s_j\ge cD^{-1/2}$ and $\sum_i c_i/s_i\le C/D$, so
$w_j^*\ge c\sqrt D$. Relative pilot accuracy and vanishing uniform mixing
preserve this bound. Set $v_L=L^{-1/16}$.
Lemma~\ref{lem:glrstable} and rounding imply that, beyond cost
$(1+\xi)L/D$ for fixed $\xi>0$, failure to certify zero requires some fresh
empirical error exceeding $v_L\sqrt D$.
Equation~\eqref{eq:checkpoint-envelope} bounds the boundary-to-cost ratio
at every later checkpoint by its value at this threshold, so the same
implication holds uniformly beyond it.
Hoeffding bounds the failure probability by
$2r\exp(-a\kappa Dv_L^2 C_s)$.
Summing geometric cost increments gives an extra expected cost bounded by
\[
 K\left(\frac L D+\frac1{\kappa Dv_L^2}\right)
 e^{-a'\kappa Lv_L^2}=o(L/D),
\]
because $\kappa Lv_L^2=L^{3/4}$ eventually. Relative checkpoint overshoot
tends to zero. Conditioning on the entire potential pilot stream is legitimate
because all main observations are fresh.

\textbf{Unsuccessful pilots and the reverse inequality.}
The contained cube of Lemma~\ref{lem:cubetail} and the checkpoint boundary
give a uniform conditional fallback bound. Put
\[
 A_L=L+\log(L+1)+\log\log(e+1/D)+\log(1/\kappa)+1,
 \qquad C_0=K_0 A_L/(\kappa D).
\]
For a sufficiently large fixed $K_0$, substitution in
\eqref{eq:checkpoint-envelope} gives $B_L(C_0)\le K_1A_L$.
Indeed $\log\log(e+C_0)$ is bounded by a constant plus the three
logarithmic terms in $A_L$; the fixed $\log\log K_0$ term grows more
slowly than $K_0$.
Rounding gives $n_j\ge (C_s-c_\Sigma)\kappa/(r c_j)$,
where $c_\Sigma=\sum_i c_i$. The shifted quotient is decreasing too:
\[
 \frac{d}{dC}\frac{B_L(C)}{C-c_\Sigma}
 =\frac{(C-c_\Sigma)B_L'(C)-B_L(C)}{(C-c_\Sigma)^2}<0,
 \qquad C>c_\Sigma.
\]
Choose $K_0$ to make every squared rectangle radius
$B_L(C_s)/(2n_j)\le b^2D/4$ at $C_0$ and, by the monotone envelope,
at all later checkpoints. The contained-cube proof then gives
$\Pp(C_{\rm main}>C_s\mid\hbox{pilot})\le2r e^{-a\kappa DC_s}$.
Summing the geometric increments yields, uniformly in the pilot outcome,
\begin{equation}\label{eq:checkpoint-fallback}
 \E(C_{\rm main}\mid\hbox{pilot})
 \le K(C_0+1/(\kappa D))\le K A_L/(\kappa D).
\end{equation}
Multiplying by $K e^{-H}$ and normalizing by $L/D$ gives $o(1)$.
The pilot and main risk budgets prove honesty. For completeness, finite
mean termination away from cancellation uses a separate fixed-instance
argument, rather than Lemma~\ref{lem:cubetail}. At any fixed strict interior
$\theta$, continuity of the fiber endpoints and
$U_\rho(P(\theta))-L_\rho(P(\theta))<2\varepsilon$ give a radius
$u_\theta>0$ such that the union of fibers over
$\|x-\theta\|_\infty\le2u_\theta$ still has width below $2\varepsilon$.
At fixed $L$, $\kappa>0$ and all main counts grow at least proportionally
to $\kappa C_s$, uniformly over frozen pilot designs. Both the checkpoint boundary and the uniform-beta
upper boundary are $O(L+\log(e+C_s))$, so eventually every rectangle
radius is at most $u_\theta$. If the empirical means are within $u_\theta$
of $\theta$, the projected rectangle then certifies an answer. Hoeffding
gives a tail at most $2r\exp(-a_\theta\kappa C_s)$ after a finite threshold,
and summing cost increments gives a finite conditional mean uniformly in
the pilot outcome. The uniform pilot rectangle has the analogous integrable
tail, so the expected pilot cost is finite as well. These constants may
depend on the fixed instance; no uniform off-cancellation efficiency is claimed.
Combining the upper bounds and letting $\xi\downarrow0$ gives the upper limit.
The moving-array lower bound in Appendix~\ref{app:uniform} needs no
logarithmic scale restriction and gives the matching lower limit.
\end{proof}

\subsection{What the joint regimes establish}
\begin{table}[tb]
\centering\small
\caption{Sufficient joint regimes under the common fixed-dimension,
fixed-cost, compact-interior and exact-cancellation assumptions.
All three conclusions concern $D_L\E C_\tau/L\to1$.
These conditions are not asserted necessary.}\label{tab:regimes}
\begin{tabular}{lll}
\toprule
Result & Pilot & Sufficient condition\\
\midrule
Theorem~\ref{thm:jointnew} & Fixed length & $L^{2/3}D_L/\log L\to\infty$\\
Theorem~\ref{thm:adaptivejoint} & Information-scaled & $\log(1/D_L)=o(L)$\\
Theorem~\ref{thm:loglogjoint} & Precision-aware & $\log(e+\log(1/g_L))=o(L)$\\
\bottomrule
\end{tabular}
\end{table}
Table~\ref{tab:regimes} distinguishes the policies and their guarantees.
Theorem~\ref{thm:loglogjoint} extends this paper's earlier sufficient regime;
it is not a new general concentration theorem. A suitable established
iterated-logarithm region can play the checkpoint region's role. The
model-specific work is the precision-aware pilot, uniform product-root
stability and expected-cost control on moving cancellation arrays.
The pilot-localized beta component is optional: removing
it leaves the theorem intact. Its numerical implementation uses integer beta
identities, outward-rounded Robbins/logarithm enclosures, and the existing
rational GLR/rectangle guards. These certify decisions within supported
integer counts, not asymptotic floating-point design optimality.

\section{Reproducible unknown-parameter experiments}
\label{sec:experimentsnew}
The computational record separates global optimization certificates from
Monte Carlo observations. Five two-ticket global certificates are supplied with executable verifiers. For example,
$\theta=(.7,.8)$, $c=(1,2)$, $\varepsilon=.22$, $\rho=.1$ has a certified
candidate $a\simeq.139421$, $w_1\simeq.553140$, and
$0.00110622158<D^*<0.00110783743$. The answer differs from both zero and
the fiber midpoint. The unknown-parameter experiments concern policies that do not receive
$\theta$, with all pilot observations charged.

\paragraph{Design, baselines, and exact returns.}
The five core instances and two stress instances are specified below.
\begin{center}\small
\begin{tabular}{lcccc}
\toprule
Instance & $\theta$ & $c$ & $\varepsilon$ & $\rho$\\
\midrule
Center & $(.5,.5)$ & $(1,4)$ & $.18$ & $.05$\\
One fair & $(.5,.8)$ & $(1,2)$ & $.18$ & $.05$\\
Robust off-center & $(.7,.8)$ & $(1,2)$ & $.22$ & $.1$\\
Exact off-center & $(.75,.8)$ & $(1,3)$ & $.24$ & $0$\\
Opposite signs & $(.2,.8)$ & $(2,1)$ & $.29$ & $.1$\\
Near boundary & $(.7,.8)$ & $(1,2)$ & $.205$ & $.1$\\
Near endpoint & $(.02,.5)$ & $(3,1)$ & $.12$ & $.02$\\
\bottomrule
\end{tabular}\end{center}
All four comparators use the same ticket access and cost accounting:
learned information allocation with the joint mixture region; equal cost
allocation with that region; a pilot Neyman allocation with that region; and
equal cost allocation with separate beta-binomial confidence sequences.
The Neyman cost weights are proportional to
$|\partial_jP(\widetilde\theta)|\sqrt{c_j\widetilde\theta_j(1-\widetilde\theta_j)}$,
with the same regularization as the information design. This is a
variance-allocation comparator, not a claimed reproduction of a published
product-estimation algorithm. Coordinate sequences each receive risk
$\delta/2$ and are projected through the exact fiber.

An initial grid at $\delta=.05$ used 80 replications per core instance and
40 per stress instance: $1{,}920$ runs in total. It discarded pilot data for
inference. That study exposed avoidable finite-confidence overhead. We then
specified the pilot-reuse variant of Corollary~\ref{cor:reuse} and evaluated
it on a separate seed range, without changing the pilot size, allocation or
boundary parameters. This follow-up has 40 replications per instance and
method, or $1{,}120$ runs. It is a fresh-seed follow-up motivated by the initial
results, not an externally preregistered study. Both complete datasets are
preserved. Equal allocation uses no pilot; learned and Neyman designs pay
for every pilot draw. The budget cap is two million queries.

Every reported return was checked with rational inequalities. Empirical
Bernoulli endpoints are handled conservatively by the fallback or coordinate
region. The sampler callback exposes only ticket responses; the policy has
no true-mean input. Runtime and optimization time are logged separately from
statistical cost. Timings were measured during parallel execution and are
not portable performance claims.

\clearpage
\begin{table}[!tb]
\centering\small
\caption{Fresh-seed pilot-reuse validation at $\delta=.05$: mean charged cost,
40 replications per cell. Learned, equal and Neyman columns use the same joint
region; Coord.\ CS uses separate coordinate mixture regions.}
\label{tab:reuse}
\begin{tabular}{lrrrr}
\toprule
Instance & Learned & Equal & Neyman & Coord. CS\\
\midrule
center & 300 & 301 & 327 & 456\\
one fair & 556 & 732 & 552 & 850\\
robust off center & 11495 & 11007 & 10890 & 15821\\
exact off center & 1043 & 874 & 972 & 1319\\
opposite sign & 3651 & 3129 & 3602 & 4859\\
near boundary & 92187 & 91756 & 98877 & 120796\\
near endpoint & 2444 & 4017 & 2254 & 2940\\
\bottomrule
\end{tabular}

\end{table}

\begin{figure}[tbp]
\centering\includegraphics[width=\linewidth]{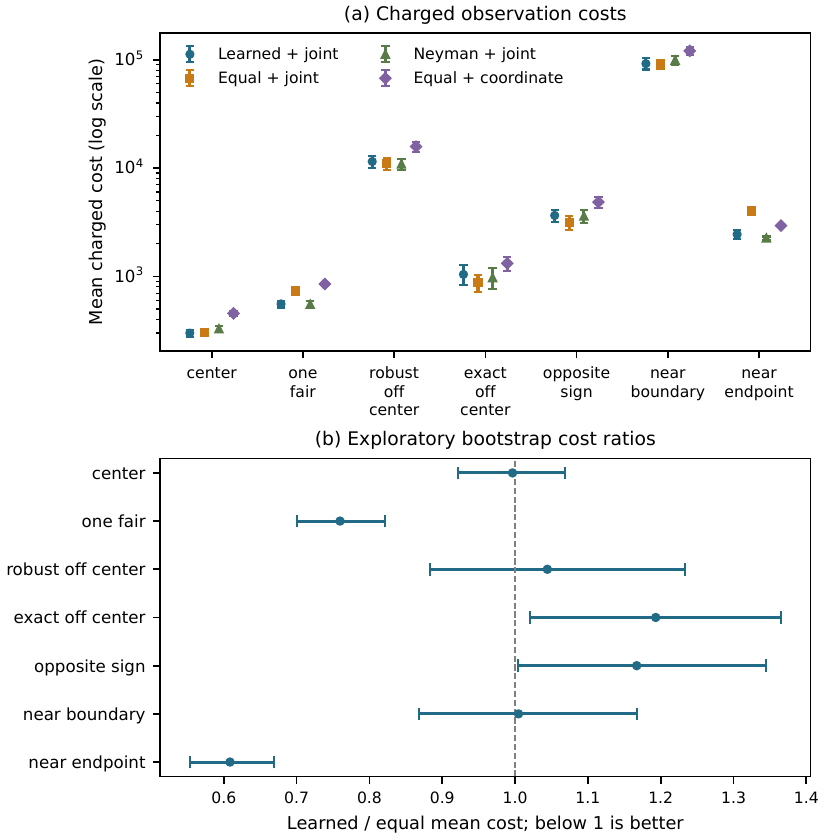}
\caption{Pilot-reuse validation. Upper panel: mean costs with $1.96$ Monte Carlo
standard errors. Lower panel: learned/equal cost ratios with exploratory 95\%
percentile bootstrap intervals, resampling replication blocks $10{,}000$
times. No multiplicity correction is applied.\label{fig:reuse}}
\end{figure}

Table~\ref{tab:reuse} and Figure~\ref{fig:reuse} show mixed performance. Learned/equal mean cost ratios
are approximately $.760$ for one fair factor and $.608$ near the Bernoulli
endpoint, but $1.193$ for the exact-sharing off-center case and $1.167$ for
the opposite-sign case. The Neyman comparator is competitive and has lower
mean cost in several cases. Reuse removes much of the full-cancellation pilot
overhead, but does not make learned information allocation uniformly best at
$\delta=.05$. These observations support the distinction between pointwise
asymptotic efficiency and ordinary-confidence performance.

\paragraph{Joint confidence--precision stress test.}
For $\theta=(.5,.8)$, unit costs, $\rho=0$, $g=1/L$, we compare the fixed
and information-scaled pilots at $L\in\{10,30,100,300,1000\}$, with 50
replications each: 500 runs. To isolate the pilot, both variants use the same
mixture boundary, regularization, fresh main data and geometric schedule.
A further 360 runs use costs $(1,3)$, $L\in\{30,100,300\}$, 30 replications
per policy, and either $\theta=(.5,.5)$ or the moving array
$\theta_L=(.5,(1+L^{-.3})/2)$. All pilot and main costs are included.

\begin{figure}[tbp]
\centering\includegraphics[width=\linewidth]{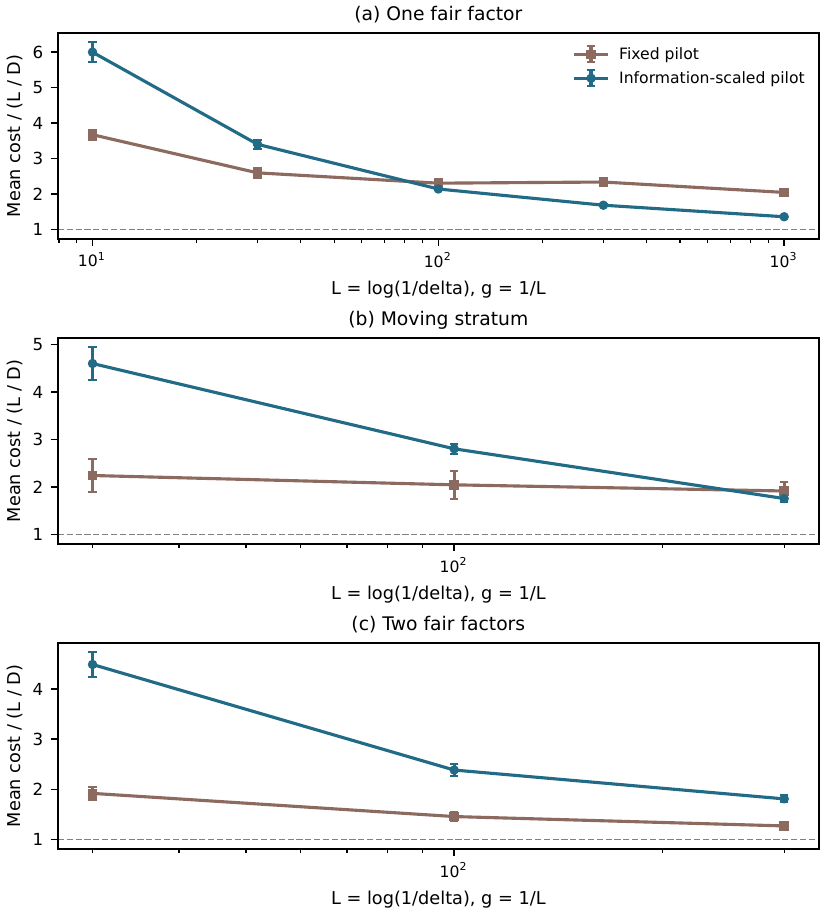}
\caption{Joint-limit diagnostics. Cost is normalized by $L/D_L$; error bars
are $1.96$ Monte Carlo standard errors. The dashed line is the asymptotic
benchmark, not a finite-$L$ attainable lower bound.\label{fig:joint}}
\end{figure}

\clearpage
Figure~\ref{fig:joint} shows that, in the one-fair case at $L=1000$, normalized means are $2.047$ for the fixed
pilot and $1.361$ for the information-scaled pilot, a mean-cost reduction of
about $33.5\%$. At $L=10$ and $30$, pilot overhead instead makes the adaptive
variant worse. In the moving-array experiment the adaptive variant improves
by $L=300$, while under full cancellation the fixed uniform-cost branch is
already well suited and is cheaper throughout the tested range. Thus the
pilot repair is justified by its regime guarantee and its targeted stress
case, not by an assertion that it should always be deployed.

All $3{,}900$ runs in the initial unknown-parameter study stopped within their caps and had no observed robust
error. Such finite experiments cannot verify risks as small as $e^{-1000}$;
the mathematical guarantee and rational return checks have distinct roles.
The historical $1{,}200$ known-instance tuned runs, 400 wrong-answer tests and
$160{,}000$ pilot-branch trials remain included and are not pooled with these
new experiments. Raw records, configurations, summaries, uncertainty
calculations and figure-generation code are provided.

\paragraph{Finite-confidence cost decomposition.}
We decomposed all 860 joint-study records using the true information $D$
and the frozen design's true zero-answer information $I_w$ for evaluation only.
With pilot cost $C_p$, main cost $C_m$ and observed boundary $\beta$,
\[
 \frac{DC}{L}=\frac{DC_p}{L}
       +\frac{D}{I_w}\frac{\beta}{L}\frac{C_m I_w}{\beta}.
\]
The last factor includes sampling fluctuation and checkpoint effects; it is
not a separately identified causal overshoot. Means of these factors should
not be multiplied as if they were independent.
At $L=1000$ on the one-fair path, the adaptive pilot's mean normalized pilot
cost is $0.160$, mean $D/I_w$ is $1.108$, mean $\beta/L$ is $1.020$,
and mean residual factor is $1.063$, with mean total $1.361$.
The fixed pilot has negligible pilot cost but mean $D/I_w=1.938$.
All observed final checkpoint cost ratios are approximately $1.15$ in these
cells. Complete decompositions, medians and 90th percentiles are supplied.

The implemented regularization remains $0.2$ through $L=390625$,
and checkpoint growth remains $1.15$ through $L\simeq1975$.
The ordinary catalog has 13 weights through $L=815730721$.
Thus the reported curves illustrate a repair but do not validate convergence
in the eventual vanishing-regularization regime. At one-fair cancellation,
$0.2$ mixing alone gives limiting fair-coordinate weight $0.9$ and leading
cost inflation $1/0.9$, before the other overheads.
A deterministic calibration compares mixing levels $0$, $0.05$ and $0.2$
at the four larger tested confidences, without new observations.

After correcting budget handling and refining numerical return bounds,
we ran 120 additional fresh-seed simulations: the one-fair, exact-sharing
off-center and near-endpoint cases, all four methods, ten replications each.
These are exploratory software and behavior checks, not a powered superiority
study. All stopped with certified returns and no observed robust error.
The supplied summaries include standard errors, medians and 90th percentiles.
Zero failures in ten replications cannot verify a five-percent risk guarantee.
The earlier fresh-seed follow-up and extensions were research development,
not externally preregistered confirmatory studies. Shared seed identifiers
form replication blocks; different binomial batch sizes do not ensure
identical observation prefixes. We make no empirical dominance claim over
the general exploration algorithms in the cited literature.

\paragraph{Precision-aware and pooled certification experiments.}
A further exploratory campaign contains 960 algorithm runs: eight cases,
four methods and 30 replications per cell. The primary comparison changes
only the main boundary of the previous adaptive-pilot policy. The remaining
methods pair the calibrated pilot with either checkpoint-only or hybrid
certification. A 480-run follow-up adds pooled rescue to both hybrid policies
using the same seeds for pathwise comparison. The follow-up was defined during
research development; it is not an independent holdout or preregistered study.
All configurations, code hashes, seeds and raw records are supplied.

For the one-fair rows in Table~\ref{tab:v09}, $\theta=(.5,.8)$ and $c=(1,1)$.
The two-fair row has $\theta=(.5,.5)$; the near-stratum row has $(.5,.516)$;
the unequal-cost row has $c=(1,4)$ and $(.5,.8)$.
The nonzero control has $\theta=(.8,.8)$, $c=(1,1)$ and
$\varepsilon=.25$. All cases in that campaign use $\rho=0$.
The total query cap is $2^{44}$ and includes all pilot observations.

\begin{table}[tb]
\centering\small
\caption{Mean charged observation cost relative to the adaptive-pilot cold baseline (C), 30 replications per cell. H changes only the main boundary; HP adds pooled rescue; CP uses the calibrated pilot, hybrid boundary and pooled rescue. Lower is better. All methods and cases are retained.}\label{tab:v09}
\begin{tabular}{lrrrrr}
\toprule
Case & $g$ & $L$ & H/C & HP/C & CP/C \\
\midrule
One fair, coarse & $5\cdot10^{-2}$ & $10$ & 0.892 & 0.642 & 0.708 \\
One fair, small & $10^{-3}$ & $10$ & 0.811 & 0.673 & 0.953 \\
One fair, tiny & $10^{-4}$ & $10$ & 0.779 & 0.651 & 0.805 \\
One fair, higher L & $10^{-4}$ & $30$ & 0.802 & 0.700 & 0.914 \\
Two fair & $10^{-3}$ & $10$ & 0.854 & 0.497 & 0.497 \\
Near stratum & $10^{-3}$ & $10$ & 0.856 & 0.594 & 0.614 \\
Unequal costs & $10^{-3}$ & $10$ & 0.912 & 0.841 & 1.231 \\
Nonzero control & $1/4$ & $\log20$ & 1.000 & 1.000 & 1.000 \\
\bottomrule
\end{tabular}
\end{table}

Changing only the boundary reduces mean charged cost by 8.8--22.1\% on the
seven cancellation cases. Adding pooled rescue to that same pilot gives
15.9--50.3\% reductions relative to the original cold policy; the nonzero
control is unchanged. Pooled rescue never increases cost in any of the 480
coupled comparisons against its own hybrid base; 399 are strict improvements.
The calibrated hybrid with rescue still costs 23.1\% more than the original
cold policy in the unequal-cost case. Its wider asymptotic guarantee is
therefore not a universal finite-sample performance claim.

All 1,440 additional runs return certified answers and have no observed
model error. This does not verify the nominal risks. The accompanying summaries
report medians, 90th percentiles and 10,000-replicate bootstrap intervals for
mean-cost ratios, resampling replication blocks within each case. These
intervals are descriptive and not adjusted for multiple comparisons.
The simulations concern supported finite counts and do not test the
exponential or double-exponential arrays in Theorem~\ref{thm:loglogjoint}.
The pooled channel's guarantee is for observation cost; extra certificate
checks can increase computation, and recorded parallel wall times are not
portable speed benchmarks.

\paragraph{Boundary and reuse ablation with common potential streams.}
To separate the effects hidden in the preceding hybrid comparison, a further
exploratory study fixes the adaptive pilot and compares cold, localized-only,
checkpoint-only and hybrid main boundaries, each with and without pooled
rescue. There are six cases and 20 replications per method, or 960 runs.
Every main region receives risk $\delta/2$; the hybrid divides that share
three ways. Pooled rescue continues the pilot process at its original
$\delta/2$ risk. The query cap is $2^{28}$, with every pilot draw charged.

The first three cases have $\varepsilon=.01$, $\rho=0$, means
$(.5,.8)$, $(.5,.5)$ and $(.5,.8)$, respectively, and costs $(1,1)$,
$(1,1)$ and $(1,4)$. The positive-radius case keeps the first means and
unit costs, with $\rho=.1$ and $\varepsilon=.109$; its zero-answer threshold
is again $g=.01$, so it is a consistency check at an equivalent threshold.
The near-nonzero case has means $(.503,.8)$, unit costs,
$\varepsilon=.01$ and $\rho=0$, giving product $.0036$ rather than exact
cancellation. The last case has means $(.7,.8)$, unit costs,
$\varepsilon=.18$ and $\rho=.05$; its product is $.24$ and zero is not
an acceptable answer. All cases use $L=10$.

Each ticket's potential Bernoulli sequence is generated in fixed blocks
keyed by seed, ticket and block index. Thus all methods within a replication
see identical prefixes even when batch sizes differ. This is stronger than
using a common seed with different binomial batch sizes. The full protocol,
source hashes, final prefixes, raw records and summaries are included.
Recovery reruns use the same deterministic study; they add no independent
observations. This is not an externally preregistered experiment.

\begin{table}[tb]
\centering\scriptsize
\caption{Boundary and reuse ablation: ratios of mean charged cost, 20 common-stream replications per cell. C: cold; L: localized; K: checkpoint; H: hybrid; P: pooled rescue. Values above one indicate higher cost.}\label{tab:v10}
\begin{tabular}{lrrrrrrr}
\toprule
Case & L/C & K/C & H/C & CP/C & LP/CP & KP/CP & HP/CP\\
\midrule
One fair & 0.829 & 1.101 & 0.850 & 0.671 & 1.000 & 1.000 & 1.000\\
Two fair & 0.894 & 1.158 & 0.912 & 0.583 & 1.000 & 1.000 & 1.000\\
Unequal costs & 0.918 & 1.060 & 0.931 & 0.832 & 1.000 & 1.000 & 1.000\\
Positive radius & 0.829 & 1.101 & 0.850 & 0.671 & 1.000 & 1.000 & 1.000\\
Near nonzero & 0.915 & 1.049 & 0.924 & 0.824 & 1.000 & 1.000 & 1.000\\
Nonzero, positive radius & 1.000 & 1.000 & 1.000 & 1.000 & 1.000 & 1.000 & 1.000\\
\bottomrule
\end{tabular}
\end{table}

Table~\ref{tab:v10} shows that localized-only certification saves
8.2--17.1\% against fresh-main cold certification in the first five cases.
Checkpoint-only instead costs 4.9--15.8\% more at this moderate confidence.
Hybrid costs are slightly above localized-only, consistent with the extra
risk split. These comparisons concern the same pilot and allocation.

The stronger reuse comparator changes the interpretation: all four pooled
variants have identical per-ticket counts and charged costs in every
replication. All 480 pooled runs return through the pooled channel, 320
before or at the corresponding base pilot's completion. Cold-pooled alone
saves 16.8--41.7\% against fresh-main cold in the first five cases; the final
case is unchanged. Hence this grid provides no incremental observation-cost
benefit for localization or hybridization beyond cold-pooled rescue.
The separate mechanisms should not receive duplicate credit for the same
savings. This negative result limits the empirical claim; it does not
invalidate the calibrated policy's asymptotic theorem.

All 960 runs stop with rationally certified answers and no observed exact
model error. The supplied summaries include means, medians, 90th percentiles
and 10,000-draw paired bootstrap intervals for cost ratios. Intervals are
descriptive, without multiplicity adjustment; cases and methods are coupled.
The study adds finite positive-radius and nonzero checks, not general
operational validation or evidence at the theorem's extreme precision scales.

\section{Discussion and conclusions}
The explicit information geometry of the partially identified product yields
computable optimal allocations, rationally verifiable two-ticket projections,
and an information-scaled pilot with a quantified joint-limit guarantee.
The counterexample to a fixed pilot and Theorem~\ref{thm:adaptivejoint}
show why learning the design must be analyzed relative to vanishing information,
rather than through fixed-precision continuity alone. Theorem~\ref{thm:loglogjoint}
extends the sufficient regime using a precision-aware pilot and a smaller
checkpoint overhead. Proposition~\ref{prop:pooledrescue} separately protects
observation cost on each coupled path, without asserting a computation-time gain.

General continuous-answer exploration is established by earlier work, including
the published AISTATS treatment of \citet{poiani2026}. The nonzero-product
fixed-width formulation of \citet{zheng1998} and the fixed-budget zero-factor
allocation of \citet{benkamra2022} are distinct precedents. The contribution
here is the explicit partial-fiber calculation, certified computation and
stated joint-limit adaptation result within this reference model.

The unknown-parameter experiments demonstrate both useful improvements and
finite-confidence overhead. Learned allocation does not uniformly outperform
equal or variance-based allocation. The information-scaled pilot is effective
in the difficult one-fair joint-limit case at high confidence, while a fixed
pilot is cheaper in several tested moderate-confidence or full-cancellation
settings. The common-stream ablation also finds no incremental hybrid gain
beyond cold-pooled rescue on its tested grid. Selection of a method should
therefore reflect its intended regime.

Efficient general-dimensional certification, a sharp joint-limit frontier,
unknown sharing radii and dependent reference answers remain open extensions.
The scalar target does not control full trajectory or perception drift.
Operational use requires a validated query protocol and an appropriate
application-level fidelity metric. Appendix~\ref{app:discovery}
clarifies discovery cost only under its explicitly restricted access model.

\begin{appendix}
\section{What passive access actually costs}
\label{sec:protocolnew}
\label{app:discovery}
A logarithmic information rate does not alone determine the number of expert
interventions when finding an auditable history is costly. Consider the
following fully specified one-ticket protocol. Each inspected inexpensive
log independently exposes a designated history with known probability
$q\in(0,1]$. Inspection costs $d\ge0$. Only when that history is exposed
may one pay $c>0$ for a fresh $\operatorname{Bern}(\theta)$ reference response.
Other log contents carry no information about $\theta$.
A reset protocol instead provides the designated ticket directly at cost $c$.
Both protocols have the same unobserved gate and accuracy task.

\begin{proposition}[The exact price of discovery]\label{prop:discovery}
Let $d_1^*(\theta)$ denote the one-ticket characteristic information per
response, with unit query cost, from \eqref{eq:generalanswergame}. At a strict
interior instance the optimal high-confidence constants are
\[
 \frac{c+d/q}{d_1^*(\theta)}\quad\text{under passive discovery},\qquad
 \frac{c}{d_1^*(\theta)}\quad\text{under direct reset}.
\]
If only expert responses are charged ($d=0$), rarity $q>0$ creates no gap
in their optimal number. If $q=0$, no data depend on $\theta$ and reliable
arbitrary-precision certification is impossible in general.
\end{proposition}
\begin{proof}
Let $M$ be inspected logs and $N$ paid expert responses. The decisions to inspect are predictable. Newly revealed exposure indicators
have conditional mean $q$ given the past. Optional sampling, first at bounded counts then by monotone
convergence, gives $\E N\le q\E M$. Thus expected cost is at least
$(c+d/q)\E N$. Exposure times and skipped logs provide independent randomness
but no parameter information, so the one-ticket lower bound still applies to
$N$. Equality is attainable in the limiting cost by using a sharp one-ticket
auditor and inspecting until each requested response becomes available:
independent geometric waiting times have mean $1/q$. The reset result is the
same theorem with direct access. At $q=0$, all transcript laws are equal
across signal parameters and the usual indistinguishable-target obstruction
applies.
\end{proof}

This specializes the two access protocols just defined; it is not a separation for all passive
logs and counterfactual continuations. In particular, claiming a $1/q$ penalty
in human interventions while treating unlimited log search as free would be
incorrect. The proposition is an elementary protocol clarification, not a new
optional-sampling principle.

\section{The earlier resolved-pilot regime and stability bounds}
\label{app:oldjoint}
This appendix concerns the fixed-length pilot specified in Section~\ref{sec:jointnew}.
\begin{theorem}[A resolved-pilot regime for simultaneous limits]
\label{thm:jointnew}
Fix $r,c,\rho$ and $a_0\in(0,1/2)$. Let
$\theta_L\in[a_0,1-a_0]^r$ satisfy $P(\theta_L)=0$, let $g_L\downarrow0$,
and let $D_L$ be the root \eqref{eq:newroot} at $(\theta_L,g_L)$.
If
\begin{equation}\label{eq:resolvedpilot}
       \frac{L^{2/3}D_L}{\log L}\longrightarrow\infty,
\end{equation}
then, at precision $\varepsilon_L=\rho+(1-\rho)g_L$,
\[
 \frac{D_L\E_{\theta_L}C_\tau}{L}\longrightarrow1.
\]
Every globally honest family has the reverse lower-limit inequality under
these arrays. At a fixed stratum with $m$ fair coordinates, a sufficient
polynomial regime is $g_L=L^{-b}$ with $0<b<m/3$.
Condition \eqref{eq:resolvedpilot} is sufficient, not necessary or sharp.
\end{theorem}
\begin{proof}
All constants below can be chosen uniformly on the displayed compact cube.
For small $D$, the inverse-KL increment satisfies
\[
 v_j(D)-m_j=\sqrt{2c_j(1-m_j^2)D}\,[1+O(\sqrt D)].
\]
In particular it is bounded above and below by positive constants times
$\sqrt D$. At least one $m_j$ is zero, so $D_L\to0$. Formula
\eqref{eq:newweights} also gives $w_j^*\ge K\sqrt{D_L}$ for every $j$:
the unnormalized weights have order
$1/[\sqrt D(m_j+\sqrt D)]$, their sum is at most $K'/D$, and at least
one fair coordinate contributes order $1/D$.

\textbf{Stability at the relevant statistical scale.}
Perturb every base magnitude by $o(\sqrt D)$. The inverse-KL roots at fixed
$D$ change by $o(\sqrt D)$ uniformly, by the preceding expansion and the
bounded derivative of $\sqrt{1-m^2}$ on the compact range.
Each root magnitude is at least $K\sqrt D$, so their product changes by a
relative $o(1)$. The logarithmic elasticity of that product in $D$ is bounded
below by a positive constant: a fair coordinate alone contributes
$1/2+o(1)$, and all other contributions are nonnegative. Bracketing the
increasing root gives $\widetilde D/D\to1$.
The same expansions in $s_j=v_j(v_j-m_j)/(1-v_j^2)$ give
$\widetilde w_j/w_j^*\to1$ for every coordinate, uniformly. The perturbed
base product is $o(g)$, because at least one base factor was zero whereas its
root factor is of order $\sqrt D$; the other perturbed factors are at most
their root factors times $1+o(1)$. Thus the pilot selects the inverse-KL
branch on this event.

Choose $u_L\downarrow0$ such that
$n_0D_Lu_L^2/\log L\to\infty$, possible by \eqref{eq:resolvedpilot}.
Hoeffding shows that pilot errors exceeding $u_L\sqrt{D_L}$ have probability
smaller than any fixed inverse power of $L$. On the complement, the selected
regularized design has zero-answer information $(1-o(1))D_L$:
its unregularized weights are relatively close to $w^*$, and adding uniform
mass $\kappa=L^{-1/8}$ can reduce the information by at most the factor
$1-\kappa$.

\textbf{Main-stage stability and tails.}
For designs with $w_j\ge K\sqrt D$, any alternative with weighted KL at
most $2D$ has coordinate KL at most $K'\sqrt D$ and hence remains in a
fixed interior cube. Taylor expansion of cross entropy, followed by weighted
Cauchy--Schwarz, gives for $\|p-\theta\|_\infty\le v\sqrt D$
\[
 \left|\sum_j(w_j/c_j)
       [\kl(p_j,x_j)-\kl(\theta_j,x_j)]\right|
 \le K(v+v^2)D
\]
on such alternatives. The same inequality with the roles reversed controls
empirical minimizers, so the projected information changes by $o(D)$ when
$v\to0$. This establishes relative GLR stability, rather than invoking
fixed-precision continuity in a shrinking neighborhood.

The mixed designs also satisfy $w_j\ge\kappa/r$. Choose $v_L\downarrow0$
so slowly that $\kappa L v_L^2/\log L\to\infty$ and that the stability
errors above vanish. After cost $(1+\xi)L/D_L$, nontermination on good pilots
requires a mean deviation exceeding $v_L\sqrt{D_L}$, with probability at
most $2r\exp(-K\kappa D_L v_L^2 C)$. Both admissible boundaries are
$L+O(\log(C+1))=L[1+o(1)]$ at that scale: condition
\eqref{eq:resolvedpilot} implies $\log(1/D_L)=O(\log L)$.
Summing the tail over stages gives $o(L/D_L)$, and stage overshoot is
$O(\sqrt{L/D_L})=o(L/D_L)$.

For arbitrary pilots, a uniform fallback bound is needed. A cube of radius
$b\sqrt D$ around $\theta$ has $|P(x)|<g$ for a sufficiently small fixed
$b>0$: compare its factors with $v_j(D)$, using the fair factor and the
uniform increment bounds. Taking smaller $b$ makes the inequality strict by
a fixed multiplicative factor. With weights at least $\kappa/r$, at cost
$K\kappa^{-1}D^{-1}[L+\log(1/D)+\log(1/\kappa)]$ the rectangle fits in
this cube except for an exponentially small empirical deviation. Tail
summation therefore bounds the conditional expected cost by this expression
uniformly in the pilot. Its contribution on the super-polynomially rare bad
pilots is $o(L/D)$. Pilot cost divided by $L/D$ is
$O(L^{-1/3}D)=o(1)$. Let $\xi\downarrow0$ for the upper bound.

\textbf{Uniform lower bound.}
Use the two opposite-sign inverse-KL alternatives at product $(1+t)g_L$,
with their extreme gates, for fixed $t>0$. Their per-ticket KL equals
$c_jD'_L$. Uniform interior Taylor bounds give likelihood variance per cost
at most $KD'_L$. The finite-cap likelihood argument in the proof of
Theorem~\ref{thm:fullgame}, now using only two output half-lines, has
martingale remainder $O(1/L)$ at cost $(1-\xi)L/D'_L$.
The root elasticity bound gives
$\limsup D'_L/D_L\le(1+t)^2$. Let $\xi,t\downarrow0$ to obtain the
reverse lower limit. Finally Corollary~\ref{cor:newstrata} gives
$D_L\asymp g_L^{2/m}$ for a fixed stratum, proving the polynomial example.
\end{proof}

\section{Uniform bounds and stopping-time details}
\label{app:uniform}
This appendix supplies the uniform estimates used in
Theorems~\ref{thm:fullgame}, \ref{thm:adaptivejoint}, \ref{thm:jointnew}
and \ref{thm:loglogjoint}.
All constants below may depend on fixed $r$, positive $c$, and an interior
compactness constant, but not on the array index, the pilot outcome or the
identity of the fair coordinate. Algorithms use parameter-independent
randomization, predictable ticket choices, stopping times, and measurable
outputs. A proposed ticket that would exceed a finite cost cap is suppressed.

The dependencies are as follows. Lemma~\ref{lem:cap} supplies lower bounds
without assuming termination at alternatives. Lemma~\ref{lem:rootstable}
transfers pilot error into relative design accuracy, and
Lemma~\ref{lem:glrstable} transfers fresh-sample error into relative evidence
accuracy. Lemma~\ref{lem:cubetail} controls the cost of a fallback even after
an inaccurate pilot. The subsequent tail sums assemble these ingredients
into the expected-cost bounds. Each local lemma states its own domain;
fixed-instance termination outside cancellation is treated separately in
the proof of Theorem~\ref{thm:loglogjoint}.

\begin{lemma}[Finite-cap likelihood bound]\label{lem:cap}
Let $E_h$, $1\le h\le M$, cover the finite-return event of interest at a
base parameter. Suppose each event is an error for a fixed alternative
$x_h$, has probability at most $\delta$ there, and has likelihood drift at
most $AC_\tau$ under the base law. If the centered log-likelihood processes
have conditional variance at most $V$ per cost unit, then
\[
 \Pp(C_\tau\le T,\ \cup_hE_h)
 \le M\delta e^{AT+b}+MVT/b^2,\qquad b>0.
\]
\end{lemma}
\begin{proof}
Run the procedure with the cap just specified. There are at most
$\lfloor T/\min c_j\rfloor$ observations. For each alternative, write
$Z_h=\sum_jN_j\kl(\theta_j,x_{hj})+M_h$ for the log likelihood of the
base relative to that alternative. The capped martingale has second moment
at most $VT$; Doob bounds the probability its maximum exceeds $b$ by
$VT/b^2$. On $E_h\cap\{C_\tau\le T,M_h\le b\}$ the likelihood is at most
$e^{AT+b}$. Integrating this likelihood against the alternative law gives
$\delta e^{AT+b}$. Sum over $h$. For the finite cover in
Theorem~\ref{thm:fullgame}, assign each output and realized cost vector to
the first covering neighborhood in a fixed finite ordering. These are
measurable events. Assign a zero-cost return an arbitrary cost vector.
Add the base-fiber error probability $2\delta$ for outputs outside
$A_\theta$. Thus no conditioning on an exact continuously distributed
answer, or termination at an alternative, is required.
\end{proof}

\begin{lemma}[Uniform inverse-root and weight stability]\label{lem:rootstable}
Fix $0<a_0<1/2$, $r\ge1$ and $c_j>0$. There exist $g_0,u_0,K>0$,
depending only on $(a_0,r,c)$, with the following property. For every
$0<g<g_0$ and $m\in[0,1-2a_0]^r$ having at least one zero coordinate,
let $D>0$ solve $\prod_jv_j(m_j,D)=g$. For every
$\widetilde m\in[0,1-a_0]^r$ and $0\le u\le u_0$, if
$\|\widetilde m-m\|_\infty\le u\sqrt D$, then the perturbed
product is below $g$, its root $\widetilde D$ exists, and
\[
 |\log(\widetilde D/D)|+
 \max_j|\log(\widetilde w_j/w_j^*)|\le Ku.
\]
The assertion also holds if the perturbation hypothesis uses
$u\sqrt{\widetilde D}$ and the perturbed root is known to exist.
\end{lemma}
\begin{proof}
Enlarge the compact magnitude range slightly. Uniform Taylor expansion of
the integral $k(m,m+q)=\int_0^q t/[1-(m+t)^2]\,dt$ gives
\[
 q=\sqrt{2c_j(1-m^2)d}\,[1+O(\sqrt d)],\qquad
 K^{-1}\sqrt d\le q\le K\sqrt d.
\]
Implicit differentiation gives
\[
 \partial_m v=
 \frac{1-v^2}{v-m}\int_m^v\frac{dt}{1-t^2}=1+O(\sqrt d),
 \qquad \partial_d v=\frac{c_j(1-v^2)}{v-m}.
\]
Consequently $|\partial_m\log v|\le K/\sqrt d$.
Writing $s=v(v-m)/(1-v^2)$, the same identities give
$|\partial_m\log s|\le K/\sqrt d$ and
$|\partial_{\log d}\log s|\le K$.
These follow also directly from $\partial_m(v-m)=O(\sqrt d)$
and $v\ge v-m\ge K^{-1}\sqrt d$.

For $H_m(d)=\prod_jv_j(m_j,d)$, the fair coordinate gives
\[
 \partial_{\log d}\log H_m(d)\ge\tfrac12-O(d),
 \qquad \partial_{\log d}\log H_m(d)\le K.
\]
The lower inequality uses the exact fair formula
$v^2=1-e^{-2c_jd}$; all other terms are nonnegative.
On $d\in[D/2,2D]$, integration of the magnitude derivatives bounds
$|\log H_{\widetilde m}(d)-\log H_m(d)|$ by $Ku$.
A sufficiently small fixed $u_0$ first brackets the perturbed root in
$[D/2,2D]$; a second bracketing gives the stated $Ku$ log-root error.
The perturbed base product is below $g$: its fair factor is at most
$u\sqrt D$, whereas the corresponding root factor is at least
$K^{-1}\sqrt D$; the other perturbed factors are at most
$(1+Ku)v_j(m_j,D)$. Reduce $u_0$ if necessary.
Integrating the derivatives of $\log s_j$ along the magnitude and root
changes, then normalizing $c_j/s_j$, proves the relative weight bound.
This argument treats $m_j/\sqrt D$ tending to zero, a positive value,
or infinity in the same inequalities.

For the last assertion work first at $d=\widetilde D$.
The product comparison is $e^{\pm Ku}$, and the true vector still has
logarithmic elasticity bounded below. Also
$g=H_{\widetilde m}(\widetilde D)\ge K^{-1}\widetilde D^{r/2}$
forces $\widetilde D\to0$. Bracket $D$ about $\widetilde D$ and apply
the preceding derivative bounds. Clipped pilot magnitudes are in the
enlarged compact range on the good event used in the main proof.
\end{proof}

\begin{lemma}[Relative GLR stability]\label{lem:glrstable}
Fix $0<a_0<1/2$, $r\ge1$, $c_j>0$ and $b_0>0$.
For $\theta\in[a_0,1-a_0]^r$ with at least one fair
coordinate, let $m_j=|2\theta_j-1|$ and let $D$ be as above.
There exist $g_0,v_0,K>0$, depending only on $(a_0,r,c,b_0)$,
such that the following holds for every $0<g<g_0$.
Let $w\in\{z\in[0,1]^r:\sum_jz_j=1\}$ satisfy $w_j\ge b_0\sqrt D$.
Define $J_w(p)=\inf_{|P(x)|\ge g}\sum_j(w_j/c_j)\kl(p_j,x_j)$.
If $J_w(\theta)\le D$, $p\in[0,1]^r$, and
$\|p-\theta\|_\infty\le v\sqrt D$ for $0\le v\le v_0$, then
\[
 |J_w(p)-J_w(\theta)|\le K(v+v^2)D.
\]
The result applies to relatively accurate optimal weights mixed with
vanishing uniform mass.
\end{lemma}
\begin{proof}
Put $a_j=w_j/c_j$. If either weighted KL objective is at most $2D$,
its individual coordinate KLs are at most $K\sqrt D$. The base and empirical
means are interior, so such alternatives lie in a common compact cube.
Taylor expansion in the first argument on that cube gives
\[
 \left|\sum_j a_j[\kl(p_j,x_j)-\kl(\theta_j,x_j)]\right|
 \le K\sum_j a_j|p_j-\theta_j||x_j-\theta_j|
       +K\sum_j a_j|p_j-\theta_j|^2.
\]
Pinsker and weighted Cauchy--Schwarz bound this by $K(v+v^2)D$.
If the empirical objective is the one at most $2D$, use
$|x-\theta|^2\le2|x-p|^2+2|p-\theta|^2$ for the same conclusion.
A true near-minimizer therefore gives $J_w(p)\le D+K(v+v^2)D<2D$;
empirical near-minimizers now satisfy the reverse inequality.
Let their approximation error decrease to zero. No compactness claim is
being made for all alternatives, only for the relevant near-minimizers.
The saddle alternative of Theorem~\ref{thm:generalgame} gives
$J_w(\theta)\le D$ for every weight vector.
\end{proof}

\begin{lemma}[Contained cube and uniform fallback]\label{lem:cubetail}
Fix $0<a_0<1/2$, $r\ge1$, positive costs $c_j$, and $0\le\rho<1/2$.
There exist $g_0,b,K>0$, depending only on these fixed quantities,
such that the following holds. Let $\theta\in[a_0,1-a_0]^r$ satisfy
$P(\theta)=0$, let $0<g<g_0$, and let $D$ be its inverse-KL product root.
Choose $g_0$ small enough that $D\le1$ uniformly. Set
$\varepsilon=\rho+(1-\rho)g$, $L\ge\log4$, and $0<\kappa\le1$.
Then
$\|x-\theta\|_\infty\le b\sqrt D$ implies $|P(x)|\le g/4$.
For every frozen cost design with $w_j\ge\kappa/r$, consider fresh counts
$n_j(s)=\lceil\alpha_jN_s\rceil$, where
$\alpha_j=(w_j/c_j)/\sum_i(w_i/c_i)$. The increasing
integer virtual totals satisfy $1\le N_0\le16$, tend to infinity,
and obey $N_{s+1}/N_s\le4$;
this includes the geometric and quadratic schedules used in the paper.
The projected rectangle from the uniform-beta region at risk $e^{-L}/2$
has conditional expected main cost at most
\[
 K\kappa^{-1}D^{-1}
 [L+\log(1/D)+\log(1/\kappa)+1].
\]
\end{lemma}
\begin{proof}
Each perturbed magnitude is at most $m_j+2b\sqrt D\le(1+Kb)v_j(D)$.
For a fair coordinate its ratio to $v_j(D)$ is at most $Kb$.
Thus the product ratio is at most $Kb(1+Kb)^{r-1}$; choose $b$ small.
The projected fiber union for $|P|\le g/4$ has width at most
$2\rho+(1-\rho)g/2<2\varepsilon$.

At total cost $C$, integer rounding gives
$n_j\ge(C-\sum c_i)w_j/c_j$. For
$C\ge C_0=K\kappa^{-1}D^{-1}
[L+\log(1/D)+\log(1/\kappa)+1]$,
the bound $\beta\le L+\log2+\sum_j\log(n_j+1)$ ensures all radii are
at most $b\sqrt D/2$. This follows by substituting $C_0$ in
$\log C$ and using that $(L+\log C)/C$ is eventually decreasing.
If every empirical error is at most $b\sqrt D/2$, the rectangle fits in the
contained cube. Hoeffding therefore bounds failure to stop after $C_0$ by
$2r e^{-a\kappa DC}$. Summing stage increments, or integrating this decreasing
tail, costs at most $K(C_0+1/(\kappa D))$. Bounded virtual-count ratios
give uniformly bounded cost ratios, so stage overshoot changes only $K$.
All statements are conditional on an arbitrary stopped pilot and use fresh
main observations, so the bound is uniform in its outcome. The restriction
$g<g_0$ is essential to these local estimates; no bound with the displayed
$\log(1/D)$ term is asserted for arbitrary large $D$.
\end{proof}

\paragraph{Pilot duration in Theorem~\ref{thm:adaptivejoint}.}
The all-round event in that proof gives
$\|\widetilde m-m\|_\infty\le K L^{-1/8}\sqrt{\widehat D}$
at every trigger. Lemma~\ref{lem:rootstable} proves relative accuracy there.
For a deterministic round with $n_kD\ge K\Lambda_{L,k}$, an error at most
$b\sqrt D$ implies $|P(\widetilde\theta)|<g/2$ by the cube comparison
and $\widehat D\ge D/2$ by root stability, after reducing $b$.
Increase $K$ so this implies the trigger. For all large $L$,
$\Lambda_{L,k+1}/\Lambda_{L,k}<2$ for every $k$, so the implication
holds at every round after $k_*$. Thus
\[
 \E n_{\rm pilot}
 \le 2^{k_*}+\sum_{k\ge k_*}2^{k+1}\,2r e^{-a2^kD}
 =O(2^{k_*}+D^{-1})=O(\sqrt L/D).
\]
The last sum is bounded by a constant times the integral of $e^{-aDx}$.
The definition of $k_*$ and $\log(1/D)=o(L)$ give
$k_*=O(L)$, $\log^2(k_*+2)=o(\sqrt L)$ and
$2^{k_*}=O(\sqrt L/D)$, including the factor-two doubling overshoot.
An earlier fallback can only shorten the pilot.

\paragraph{Main-stage integration and unsuccessful pilots.}
On the good pilot event, Lemma~\ref{lem:rootstable} and uniform mixing give
$J_w(\theta)\ge(1-o(1))D$ and $w_j\ge K^{-1}\sqrt D$.
Rounding the counts changes relative weights by $O(1/(C\sqrt D))=o(1)$
at $C\asymp L/D$. Lemma~\ref{lem:glrstable} applies.
For any fixed $\xi>0$, choose $v_L=L^{-1/16}$. At
$C_1=(1+\xi)L/D$, the boundary is $L+o(L)$ and
nontermination requires an empirical error exceeding $v_L\sqrt D$.
Its tail is at most $2r e^{-a\kappa Dv_L^2 C}$.
The cost beyond $C_1$, apart from relative stage overshoot tending to zero,
is bounded by
\[
 K\left(C_1+\frac1{\kappa Dv_L^2}\right)
 \exp\{-a'\kappa Dv_L^2 C_1\}=o(L/D),
\]
since $\kappa L v_L^2=L^{3/4}$ eventually.
By Lemma~\ref{lem:cubetail}, the bad-pilot contribution divided by $L/D$
is at most
\[
 K e^{-L^{1/4}}\kappa^{-1}
 \left[1+\frac{\log(1/D)+\log(1/\kappa)+1}{L}\right]=o(1).
\]
Realizing complete potential pilot streams independently of main streams
justifies this multiplication even though the good event includes unused
pilot observations. Add the pilot expectation and let $\xi\downarrow0$.

\paragraph{Moving-array lower bound.}
Choose two inverse-KL alternatives at product magnitude $(1+t)g$,
with opposite signs obtained by flipping a truly fair coordinate.
Every coordinate KL is $c_jD'$. Their extreme targets are
$\pm[\rho+(1-\rho)(1+t)g]$, beyond $\pm\varepsilon$.
The negative target makes every nonnegative output erroneous; the positive
target makes every negative output erroneous. These two events cover all
finite returns. Interior Taylor bounds give likelihood variance per cost
at most $KD'$. Lemma~\ref{lem:cap}, with
$A=D'$, $V=KD'$, $T=(1-\xi)L/D'$ and $b=\xi L/2$, has remainder $O(1/L)$.
Hence the expected cost is at least $(1-o(1))(1-\xi)L/D'$.
The fair-coordinate elasticity gives
$\limsup D'/D\le(1+t)^2$ by integrating $d\log H/d\log d\ge1/2-o(1)$.
Let first $L\to\infty$, then $\xi,t\downarrow0$.
This lower bound does not require the older pilot-resolution condition.

\paragraph{Pooled pilot data.}
For Corollary~\ref{cor:reuse}, apply Ville to
$\exp\{4u(S_n-n\theta)-2nu^2\}$ and its lower-tail counterpart.
Crossing a deviation $u$ at any $n\ge m$ forces a value at least
$e^{2mu^2}$, proving the stated uniform prefix bound.
For an arbitrary fixed pilot, write the pooled centered mean as
$n_0(\widehat\theta_0-\theta)/n+(n-n_0)(\widehat\theta_{\rm fresh}-\theta)/n$.
Once $n\ge2n_0/u$, the first term is at most $u/2$ in magnitude.
For fixed sufficiently small $u$ at a strict instance, positive regularization
and Hoeffding therefore give conditional fallback cost
$K\kappa^{-1}(L+n_0/u)=O(\kappa^{-1}(L+n_0))$.
On good pilots use the unconditional prefix bound with
$n\ge K^{-1}\kappa C$, rather than conditioning on reused observations.
Multiplying the conditional fallback bound by the exponentially small
bad-pilot probability proves the claimed expectation result.

\paragraph{Allocations along moving cancellation strata.}
In the notation following Corollary~\ref{cor:newstrata}, put
$H=\{j:\beta_j\ge\gamma\}$ and
\[
 A_j=\begin{cases}
 (2d)^{-1},&\beta_j>\gamma,\\
 c_j/\{(a_j+\sqrt{2c_jd})\sqrt{2c_jd}\},&\beta_j=\gamma.
 \end{cases}
\]
Then $w_j\to A_j/\sum_{i\in H}A_i$ on $H$.
For $\beta_j<\gamma$, $w_j\sim
 B_jg^{\gamma-\beta_j}/\sum_{i\in H}A_i$, where
$B_j=\sqrt{c_j/(2d)}/a_j$ if $\beta_j>0$, and
$B_j=\sqrt{c_j(1-a_j^2)/(2d)}/a_j$ if $\beta_j=0$.
These formulas follow directly by substituting the uniform increment
expansion in $c_j/s_j$ and normalizing. They display the variance factor
needed for a nonvanishing biased magnitude.

\end{appendix}

\begin{acks}[Use of generative AI]
OpenAI ChatGPT/Codex assisted with mathematical exploration and draft proofs,
software implementation and debugging, experiment execution, literature
searches, manuscript drafting and journal formatting. Its role extended beyond
language editing. The authors are responsible for the accuracy, originality,
attribution and integrity of the final scientific claims.
\end{acks}
\begin{funding}
The authors received no funding for this research.
\end{funding}
\section*{Competing interests}
The authors declare that they have no competing interests.
\section*{Data and code availability}
All data are synthetic. Executable simulation and analysis code, complete raw
records for 6,420 algorithm runs (3,900 historical, 120 earlier revision, 1,440 localization runs, and 960 boundary/reuse ablation runs), fixed seeds, configurations, figures and
verification scripts are supplied in Supplement A. No human participants or
personal data were used in these experiments.
\begin{supplement}
\stitle{Supplement A: reproducibility code and synthetic results}
\sdescription{Self-contained source code, raw simulation records, configurations,
summary and figure-generation scripts, tests, exact certificate verifiers and
reproduction instructions. The archive is supplied as an accompanying file,
not exclusively as an external link.}
\end{supplement}
\bibliographystyle{imsart-nameyear}
\bibliography{references}
\end{document}